\documentclass[mnsc,nonblindrev]{informs3mod} %

\OneAndAHalfSpacedXI %

\usepackage[colorlinks=true]{hyperref}
\usepackage{bm}
\usepackage{verbatim}
\usepackage{tikz}
\usepackage{pgfplots,pgfplotstable}
\pgfplotsset{compat=1.3}
\usetikzlibrary{shapes.geometric, arrows, patterns}
\tikzstyle{box} = [rectangle, rounded corners, minimum width=2cm, minimum height=1cm, text centered, draw=black, align=center]
\tikzstyle{arrow} = [thick,->,>=stealth]
\tikzstyle{lrarrow} = [thick,<->,>=stealth]
\tikzstyle{darrow} = [thick,dotted,->,>=stealth]
\usepackage{multicol}

\theoremstyle{definition}

\theoremstyle{remark}

\def\cbl{\color{black}}
\def\ctl{\color{black}}
\definecolor{gold}{rgb}{0.85,0.65,0}

\renewcommand{\theequation}{\mbox{\arabic{equation}}}
\def\beq{\begin{equation}}
\def\eeq{\end{equation}}

\def\fnote#1{\footnote}
\newcommand{\epr}{\hfill\hbox{\hskip 4pt \vrule width 5pt height 6pt depth 1.5pt}\vspace{0.0cm}\par}

\newcommand{\grad}{\ensuremath{\nabla}}

\def\cG{{\cal G}}

\def\cP{{\cal P}}

\newcommand{\bbE}{\mathbb{E}}

\newcommand{\bbN}{\mathbb{N}}

\newcommand{\bbP}{\mathbb{P}}
\newcommand{\bbQ}{\mathbb{Q}}
\newcommand{\bbR}{\mathbb{R}}

\newcommand{\CG}{\mathcal{G}}

\newcommand{\CP}{\mathcal{P}}

\DeclareMathOperator*{\esssup}{ess\,sup}

\DeclareMathOperator{\Proj}{Proj}

\DeclareMathOperator{\sign}{sign}

\DeclareMathOperator{\Tr}{Tr}

\DeclareMathOperator{\cl}{cl}

\DeclareMathOperator{\Conv}{Conv}

\DeclareMathOperator{\SPO}{SPO}
\DeclareMathOperator{\LS}{LS}
\DeclareMathOperator{\AD}{AD}

\DeclareMathOperator{\SPOp}{SPO+}

\DeclareMathOperator{\Bayes}{Bayes}
\DeclareMathOperator{\reg}{reg}

\DeclareMathOperator{\sym}{sym}
\DeclareMathOperator{\cont}{cont}

\usepackage{natbib}
 \bibpunct[, ]{(}{)}{,}{a}{}{,}%
 \def\bibfont{\small}%
 \def\BIBand{and}%

\TheoremsNumberedThrough     %
\ECRepeatTheorems

\EquationsNumberedThrough    %

\begin{document}

\RUNTITLE{Risk Guarantees for End-to-End Prediction and Optimization Processes}

\TITLE{Risk Guarantees for End-to-End Prediction and Optimization Processes}

\ARTICLEAUTHORS{%
\AUTHOR{Nam Ho-Nguyen}
\AFF{Discipline of Business Analytics, The University of Sydney. \EMAIL{nam.ho-nguyen@sydney.edu.au}} %
\AUTHOR{Fatma K{\i}l{\i}n\c{c}-Karzan}
\AFF{Tepper School of Business, Carnegie Mellon University. \EMAIL{fkilinc@andrew.cmu.edu}} %
} %

\ABSTRACT{%
Prediction models are often employed in estimating parameters of optimization models. Despite the fact that in an end-to-end view, the real goal is to achieve good optimization performance, the prediction performance is measured on its own. While it is usually believed that good prediction performance in estimating the parameters will result in good subsequent optimization performance, formal theoretical guarantees on this are notably lacking. In this paper, we explore conditions that allow us to explicitly describe how the prediction performance governs the optimization performance. Our weaker condition allows for an asymptotic convergence result, while our stronger condition allows for exact quantification of the optimization performance in terms of the prediction performance. In general, verification of these conditions is a non-trivial task. Nevertheless, we show that our weaker condition is equivalent to the well-known Fisher consistency concept from the learning theory literature. This then allows us to easily check our weaker condition for several loss functions. We also establish that the squared error loss function satisfies our stronger condition. Consequently, we derive the exact theoretical relationship between prediction performance measured with the squared loss, as well as a class of symmetric loss functions, and the subsequent optimization performance. In a computational study on {\cbl portfolio optimization, fractional knapsack and multiclass classification problems,} we compare the optimization performance of using of several {\cbl prediction loss functions (some that are Fisher consistent and some that are not) and demonstrate that lack of consistency of the loss function can indeed have a detrimental effect on performance.}
}%


\maketitle

\section{Introduction}\label{sec:jpo-intro}

{\cbl The optimum solutions of optimization models crucially depend on the parameters defining these models, but these parameters are hardly ever available directly.} In practice, these `true' model parameters are predicted from side information and historical data often using statistical inference or machine learning techniques. There are many techniques that quantify the performance of prediction models. Nevertheless, these techniques almost exclusively focus on achieving a good prediction performance and do not take into account the subsequent optimization task. This is despite the fact that the ultimate goal in this process is to make the best decision in the subsequent optimization problem, not necessarily to have the best generic prediction performance of the parameters. In this paper, we consider a joint  \emph{end-to-end} view of the prediction and optimization processes, and 
identify the critical properties of prediction models in terms of guaranteeing a low optimality gap in the subsequent optimization performance.

More formally, we consider {\cbl an} optimization problem of the form
\begin{equation}\label{eqn:jpo-opt-problem}
\min_x \left\{ f(x) + c^\top x :~ x \in X \right\},
\end{equation}
where $X \subset \bbR^m$ is a convex compact domain, and $f:X \to \bbR$ is a convex function (and hence continuous on the relative interior of $X$). In our setting, the linear vector $c$ is not known exactly, but instead is governed via covariates $w$. More precisely, we suppose the covariates $w$ belong to a given set $W \subseteq \bbR^k$, and the vectors $c$ belong to a given set $C \subseteq \bbR^m$. We assume that $(w,c) \sim \bbP$ for some unknown distribution $\bbP$ on $W \times C$, and we need to solve \eqref{eqn:jpo-opt-problem} for $c$ yet we are only given information of $w$. {\cbl Note that our setup covers the case when $c$ is still noisy even when given $w$, since the conditional distribution $\bbP[c \mid w]$ may not be a point mass, and indeed this will be the more interesting case that we will study. Note also that previous literature studied the case when the function $f = 0$, but we consider a general convex $f$ function, which is relevant in many applications; see Example~\ref{ex:jpo-port-opt}.}

While we do not know the distribution $\bbP$, we have access to {\cbl the} historical data $H_n := \{(w_i,c_i) : i \in [n]\}$, 
where %
the $(w_i,c_i)$ are realizations of independent and identically distributed (i.i.d.) random variables %
from the unknown distribution $\bbP$. We examine an end-to-end view of the following prediction and optimization processes: first, based on data $H_n$, a prediction model in the form of a function $g:W \to \bbR^m$ is built to capture the dependency of $c$ on $w$; then, when given a covariate $w$, \eqref{eqn:jpo-opt-problem} is solved with $c$ replaced by the prediction $g(w)$. This setting is commonly used amongst practitioners in decision-making domains for a variety of problems. Below, we give {\cbl three} %
particular examples, although many more exist.

\begin{example}\label{ex:jpo-shortest-paths}
Suppose we have a collection of service items (e.g., machines, vehicles) which we maintain over a certain time horizon. These items need refurbishment or replacement after a certain number of time periods. The optimal maintenance schedule can be defined as a shortest path problem over an appropriately defined network, where the `distances' are given by the maintenance costs. Note that such future costs are often obtained via forecasts, and thus are not deterministic. In this setting, $X$ is the convex hull of all paths from the starting point to the ending point in the underlying graph (each such path represents a maintenance plan), $f(x) = 0$ for all $x \in X$, and $c$ is the vector of arc distances that represent the maintenance costs. Side information (covariates) $w$ of $c$ can consist of (amongst others) seasonality, demand, supply and other economic factors.
\epr
\end{example}

\begin{example}\label{ex:jpo-port-opt}
Consider a portfolio optimization problem, where the task is to allocate wealth to $m$ different assets to maximize investment return. In the typical mean-variance formulation, the goal is to simultaneously minimize the variance of the portfolio return, while maximizing the expected return. Then, $X$ is the set of all possible asset allocations, each $x \in X$ represents an asset allocation (i.e., $x_j$ represents how much wealth to invest into asset $j$), $f(x) = \gamma x^\top \Sigma x$ is the variance of the portfolio return with $\Sigma$ being the covariance matrix of the returns between the assets, and $c = -\mu$ is the vector of (negative) returns for each assets. In many settings, $\Sigma$ is assumed {\cbl to be} stable, and $\mu$ is predicted through market factors (e.g., liquidity, value, momentum, volume) which can be considered as the side information $w$.
\epr
\end{example}

{\ctl
\begin{example}\label{ex:structured-prediction}
{\cbl Structured prediction} is a form of multiclass classification designed to predict structured objects, such as sequences or graphs, from feature data; see e.g., \citet{GohJaillet2017,osokin2017structured} and references therein. In structured prediction, given covariates $w$, a structured object $\tilde{x}$ from some output space $\tilde{X}$ is chosen as the prediction, often by solving $\min_{\tilde{x} \in \tilde{X}} \tilde{g}(\tilde{x};w)$. In this setting, $\tilde{X}$ is usually a finite combinatorial set, so $\tilde{g}(\tilde{x};w)$ is a vector {\cbl in which each coordinate corresponds to the cost of an object} %
$\tilde{x} \in \tilde{X}$. The `true' structured loss in this setting is measured between the selected $\tilde{x}$ and the `correct' $\tilde{x}^* \in \tilde{X}$, denoted {\cbl by} $L(\tilde{x},\tilde{x}^*)$. Since $\tilde{X}$ is combinatorial, $L$ can be defined as the Hamming loss, although some other structured losses are possible. Structured prediction fits into our optimization setting by taking $X$ to be a simplex whose vertices correspond to objects in $\tilde{X}$, $g(w) = \{\tilde{g}(\tilde{x};w)\}_{\tilde{x} \in \tilde{X}}$, and the Hamming loss can be cast as the optimality gap of a particular constructed cost vector $c_{\tilde{x}^*}$. 
\epr
\end{example}
}

{\ctl %

{\cbl Given a point $(w,c)\in W \times \bbR^m$ and a prediction function $g:W \to \bbR^m$, in order to} assess the quality of using $g(w)$ in place of $c$ in \eqref{eqn:jpo-opt-problem}, we define the \emph{true loss} as the optimality gap of the solution obtained with $g(w)$ on the true objective vector $c$, that is, the quality of the prediction $d = g(w)$ with respect to \eqref{eqn:jpo-opt-problem} is given by the \emph{true loss function}
\[
L(d,c) := f(x^*(d)) + c^\top x^*(d) - \min_{x \in X} \left\{ f(x) + c^\top x \right\},
\]
{\cbl where $x^*(d)\in\argmin_{x \in X} \left\{ f(x) + d^\top x \right\}$.}
Since $(w,c)$ is randomly drawn from $\bbP$, we assess the performance of a function $g:W \to \bbR^m$ in terms of the expected true loss,  i.e., the \emph{true risk}
\[R(g,\bbP) := \bbE[L(g(w),c)].\]

A na\"{i}ve attempt to minimize $R(g,\bbP)$ is to directly use empirical risk minimization (ERM) with the loss $L$ to train the prediction model $g$, i.e., given the data $H_n = \{(w_i,c_i)\}_{i \in [n]}$, we obtain a prediction function $\hat{g}$ by solving
\[ \inf_g \frac{1}{n} \sum_{i \in [n]} L(g(w_i),c_i)). \]
However, $L(d,c)$ is not convex in $d$, and thus it is not possible in general to obtain a polynomial-time approach with certified performance guarantees from minimizing the empirical risk based on the true loss function $L$. The natural remedy is to use a \emph{convex surrogate loss function} $\ell$ in place of $L$. The use of surrogate loss functions to ensure algorithmic tractability is very common in machine learning. For example, convex surrogates such as hinge loss are used instead of the non-convex true $0$-$1$ loss in classification problems.

{\cbl The question now becomes: which surrogate loss function should we use? While quite a number of surrogate loss functions have been proposed and used in this context, it is not yet well-understood how using a regression technique performs in terms of the true risk. More precisely, we define the \emph{surrogate risk} as
\[R_\ell(g,\bbP) := \bbE\left[ \ell(g(w),c) \right].\]
This paper aims to understand how a minimization scheme for the surrogate risk $R_\ell(g,\bbP)$, which is well-established and implementable in practice, can impact the true risk $R(g,\bbP)$. In other words, if we employ an established regression technique to obtain a prediction function $\hat{g}$, what can we say about the true risk $R(\hat{g},\bbP)$ of $\hat{g}$? Consequently, to fill this gap in the literature, in this paper we explore this relationship and identify important properties of surrogate loss functions $\ell$ that enable us to derive guarantees on the true risk. We make these concepts mathematically rigorous, and describe their relationship to traditional notions of statistical consistency, in Section \ref{sec:jpo-risk}. For a visual summary of our framework see Figure \ref{fig:jpo-description}.

\begin{figure}[t!bh]
	\centering
	\includegraphics[page=1,scale=1]{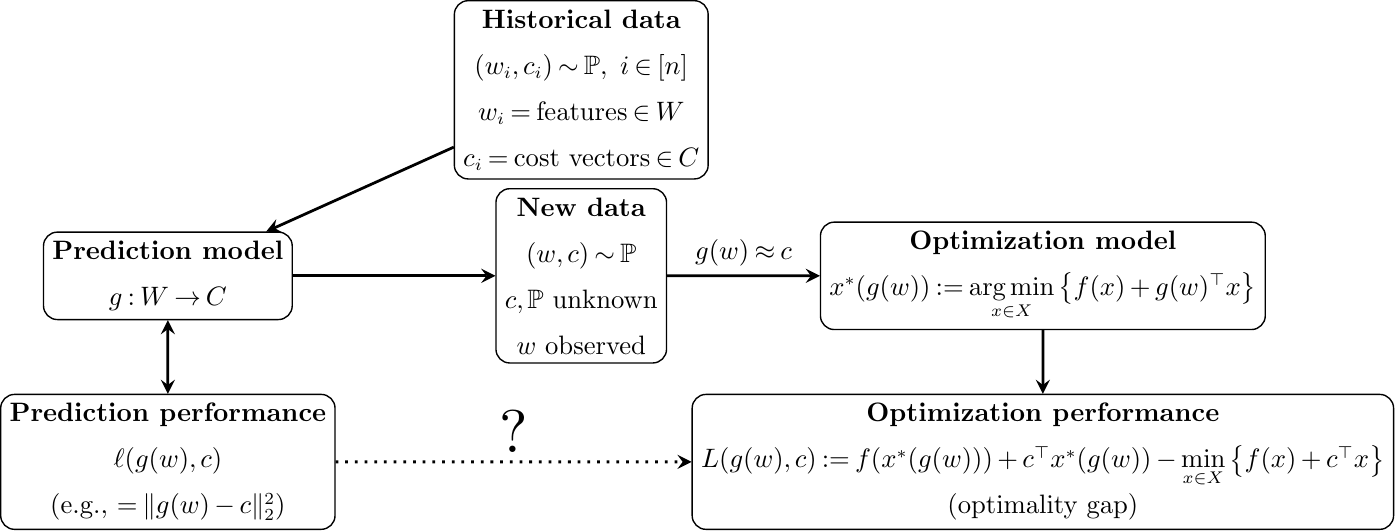}
	\caption{The end-to-end prediction and optimization framework.}\label{fig:jpo-description}
\end{figure}
}
}

\subsection{Outline and Contributions}

In this paper, we examine, from an end-to-end view, how the performance of the prediction part relates to the performance of the optimization part. In particular, we establish conditions for the existence of explicit relationships between the prediction performance, i.e., the surrogate risk, and the optimization performance, i.e., the true risk. 

{\cbl In Section \ref{sec:jpo-literature}, we review the literature related to this topic.} In Section \ref{sec:jpo-risk}, we precisely define the problem we address, and outline the challenges. 

In Section \ref{sec:jpo-FisherConsistent-surrogate-risk}, we rigorously derive {\cbl technical sufficient} conditions on the prediction loss function that allow us to asymptotically minimize the true risk by minimizing the surrogate risk for a given probability distribution $\bbP$.
These conditions are based solely on the choice of prediction loss function, rather than the class of prediction models $g$ that we wish to select from. They allow us to compare and contrast the resulting optimization performance when different prediction model training methods are used for the estimation of objective function parameters, and thus are instrumental in terms of {\cbl selecting} among such training methods. In addition, {\cbl our results in this section make the following contributions:
\begin{itemize}
	\item We show that, in the prediction and optimization context, the concept of calibration introduced by \citet{Steinwart2007}, which allows us to establish performance guarantees, is equivalent to the well-known concept of Fisher consistency. To the best of our knowledge, such a relationship was not described in the previous literature. This result provides a tool for easily checking which loss functions lead to performance guarantees, which we exploit in our examples in Section \ref{sec:jpo-FisherConsistent-surrogate-risk}.

	\item We compare several prediction methods from practice through the lens of our conditions. In Example \ref{ex:jpo-squaredLoss}, we show that the commonly used squared loss function $\ell(d,c) = \|d-c\|_2^2$ satisfies Fisher consistency. In Examples \ref{ex:jpo-spo-loss}--\ref{ex:jpo-multiclass}, we examine the $\SPOp$ loss function from \citet{ElmachtoubGrigas2017}, which is particularly relevant since it is the only convex loss function (thus far) that incorporates optimization problem information in the prediction and optimization setting. The Fisher consistency of the $\SPOp$ loss function in certain settings was previously established in \citet{ElmachtoubGrigas2017} (which we summarize in Examples~\ref{ex:jpo-spo-loss}--\ref{ex:jpo-binary-class}). Despite this, our Examples~\ref{ex:jpo-spo-one-dim}--\ref{ex:jpo-multiclass} show that the $\SPOp$ loss function is not Fisher consistent in other fairly natural settings such as multiclass classification. Furthermore, our numerical study in Section~\ref{sec:num-multiclass} highlights the importance of having Fisher consistency of a loss function over simply a property that the loss function is customized to the optimization problem.
\end{itemize}
}

Often in statistical learning, we {\cbl are given} minimal knowledge of the distribution $\bbP$. Therefore, the distribution dependent nature of our results from Section~\ref{sec:jpo-FisherConsistent-surrogate-risk} is not so desirable. In Section \ref{sec:jpo-non-asymptotic}, building on the results from \citet{Steinwart2007}, we establish conditions for distribution-independent relationships between {\cbl the} true risk and the surrogate risk. {\cbl In this section, our main contributions are as follows:}
\begin{itemize}
	\item Since checking these conditions is difficult for loss functions in general, we first focus on a tractable special case of using the squared loss function to measure prediction performance (i.e., the least squares method to train a prediction model). For the first time in the literature of joint prediction and optimization setting, using our conditions, we prove an explicit relationship between the surrogate squared risk and the true optimization risk. This then allows us to relate our true risk to a class of symmetric loss functions by exploiting existing results on {\cbl regression} from \citet{Steinwart2007}; {\cbl see Section~\ref{sec:U-calib-LS}}.
	
	\item We also study distribution-independent risk relationships for the $\SPOp$ loss function proposed in \citet{ElmachtoubGrigas2017} {\cbl in Section~\ref{sec:jpo-uniform-calibration-spo}}. {\cbl The $\SPOp$ loss has not been formally studied in this context before, and due to the importance of the $\SPOp$ loss to the prediction and optimization setting, we believe that such a study is warranted.}
\end{itemize}

In Section~\ref{sec:jpo-computational}, we carry out a computational study on {\cbl three problem classes on real and simulated data: portfolio optimization, fractional knapsack, and multiclass classification.
Our study on portfolio optimization is based on real-world data, where consistency is not known a priori. Our study with the fractional knapsack problem on simulated data allows us to chose some parameters to control the degree of non-linearity of the underlying data model, and thereby the model misspecification of certain loss functions. Lastly, we examine multiclass classification on simulated data, where the $\SPOp$ loss is provably Fisher inconsistent (see Example \ref{ex:jpo-multiclass}), but the squared loss is consistent.

Our numerical findings support our theoretical results by indicating that the conditions we identified for the loss function have a non-trivial effect on practical performance.
{\cbl On the real world instances of portfolio optimization, we observe that there is little difference between using squared loss and $\SPOp$ loss, where it is likely that both of these loss functions are consistent and there is no model misspecification. Indeed, we show in Section~\ref{sec:portfolio} that for a specific convex quadratic program with uncertain linear term and a single linear constraint and no non-negativity constraints arising in the mean-variance portfolio optimization, the $\SPOp$ loss is equivalent to the true loss $L$. Furthermore, for this problem, we show that for the class of linear predictors, the optimal least squares risk predictor is also optimal for the true risk. Because of this, we carry out our experiments on more interesting case of the portfolio optimization instances with nonnegativity constraints. 
In contrast, our experiments on multiclass classification highlight an important insight: consistency of a loss function matters as much as (if not more) whether the loss function takes into account optimization problem information. In particular, despite the fact that the $\SPOp$ loss takes into account information from the optimization problem, its inconsistency for multiclass classification  problem resulted in poor performance. On the fractional knapsack instances, we re-affirm the observation of  \citet{ElmachtoubGrigas2017} that the degree of model misspecification plays a role favoring $\SPOp$ loss over squared loss when there is no significant difference between consistency and calibration properties of the loss functions. 
}}

{\cbl We relegate all of the proofs to the appendices of the corresponding sections.}

\paragraph*{Notation.} We use of the following notation. Given a positive integer $N$, $[N] := \{1,\ldots,N\}$. Throughout, $k,m \in \bbN$ are the dimensions of the Euclidean spaces where $W,C$ live respectively, $j \in [m]$ always denotes an index for the component of a vector in $\bbR^m$, and $i \in [n]$ denotes an index for a data point $(w_i,c_i) \in H_n$. 
Given a vector $d \in \bbR^m$, {\cbl we let} $X^*(d) := \argmin_{x \in X} \left\{ f(x) + d^\top x \right\}$ %
{\cbl to be} the argmin mapping, and $x^*(d)$ denotes some selection from $X^*(d)$ {selected in a} deterministic {\cbl manner}. More precisely, $x^*:\bbR^m \to X$ is a function such that for any $d \in \bbR^m$, $x^*(d) \in X^*(d)$. Our results are agnostic to the specific choice of algorithm {\cbl picking $x^*(d) \in X^*(d)$}.

\section{Related Literature}\label{sec:jpo-literature}

Both prediction and optimization have been studied extensively on their own. In particular, the selection of the prediction function $g:W \to \bbR^m$ to minimize some measure of prediction error on the given data $H_n$ is studied extensively in statistics and machine learning, see e.g., \citet{BousquetEtAl2004}.
{\cbl
Moreover, a classical machine learning application, that is the classification problem, where a prediction model is built first from training data based on a loss function, presents a setup close to our end-to-end joint prediction and optimization view. The benchmark loss function in the context of the classification problem is the $0$-$1$ loss, but it is nonconvex. Thus, in order to get polynomial-time algorithms for training, $0$-$1$ loss is often replaced with a convex surrogate loss function. Consequently, this necessitates the study of the relationship between the surrogate loss functions and the true $0$-$1$ loss within this context. This is a topic well-studied and understood; for example, \cite{Steinwart2002a,Steinwart2002b,Lin2004,Zhang2004,Steinwart2005,BartlettJordanMcAuliffe2006} have developed a general theory for the minimization of the true $0$-$1$ risk via a surrogate risk which satisfies certain criteria. This was extended to robust regression and density estimation problems by \citet{Steinwart2007}, who builds a theory for the relationship between true and surrogate risk. Our work can be seen as a generalization of these results to optimization problems involving prediction parameters. %
In this context, our optimality gap is analogous to the $0$-$1$ loss in classification, but is much more complicated.
}

From an end-to-end point of view, the relationship between the prediction models used to obtain model parameters and the subsequent optimization performance has, to our knowledge, only been examined by a few papers. This line of work was initiated by \citet{bengio1997using} who explored the use of a financial training criterion in neural networks rather than a prediction criterion. In the context of newsvendor inventory control problem, \citet{liyanage2005practical} {\cbl show that, rather than analyzing the optimal order quantity derived for the distribution that is estimated from the data, it is better to propose a broader class of order policies and choose the optimal policy that maximizes the expected profit on the data.} 
More recently, \citet{KaoVRYan2009,ElmachtoubGrigas2017} and \citet{DontiAmosKolter2017} contributed to this line of research. These papers examined designing or using alternative loss functions in training the prediction model so as to improve the final optimization performance. \citet{KaoVRYan2009} study the specialized setting where $X = \bbR^m$, $f$ is a strongly convex quadratic, and the prediction model $g$ is restricted to be linear, and present theoretical guarantees under a particular data distribution. \citet{DontiAmosKolter2017} propose a scheme to directly differentiate the optimality gap, which gives rise to a stochastic gradient descent scheme for directly training the prediction model via the optimality gap. While superior numerical performance of this algorithmic scheme was demonstrated in \citet{DontiAmosKolter2017}, they provide no theoretical guarantees for the convergence of the risk quantities in their approach. {\cbl %
In a setting closest to ours, \citet{ElmachtoubGrigas2017} examine the true optimality gap loss, and propose a convex surrogate loss from a quantity upper bounding the true optimality gap, and suggest that this convex surrogate loss function, referred to as the $\SPOp$ loss, should be used in prediction model training.} They show {\cbl Fisher consistency (see Definition \ref{def:jpo-FisherConsistent-loss})} of their surrogate loss under certain distributional assumptions, but do not give explicit relationships on how the performance of the prediction part governs the optimization performance. {\cbl In contrast to their work on designing a new surrogate loss function, the main goal} of our paper is essentially {\cbl to close} this theoretical gap in the literature {\cbl by identifying \emph{properties} of loss functions that ensure good performance and providing explicit relationships between the performance of the prediction loss and the optimization loss for \emph{general} classes of loss functions satisfying these properties; see Sections \ref{sec:jpo-FisherConsistent-surrogate-risk} and \ref{sec:jpo-non-asymptotic}. As such the focus and the results presented in our paper are very different than the ones from  \citet{ElmachtoubGrigas2017}. Note that in certain parts of our paper, we use specific loss functions, such as the squared loss or $\SPO+$ loss to demonstrate that they possess or lack certain properties that we have identified. For this purpose, the squared loss is rather classical, and the main purpose in designing the $\SPO+$ loss in Elmachtoub and Grigas (2017) was to keep the end-to-end framework in view, and therefore it is a very natural candidate to examine.}

{\cbl As we discussed in Example~\ref{ex:structured-prediction}, this paper is also related to structured prediction.} 
\citet{osokin2017structured} provides risk relationships between surrogate methods to predict the vector $\{\tilde{g}(\tilde{x};w)\}_{\tilde{x} \in \tilde{X}}$ {\cbl (see Example~\ref{ex:structured-prediction})} and the true structured loss. Our goal in this paper is to provide results for the more general optimization setting, where there are a potentially infinite number of `objects', and when the true loss is the optimality gap.

As an alternative approach to this end-to-end view of the predict-then-optimize framework, one may wish to avoid appealing to an explicit prediction model completely, and instead use density estimation as a compelling method to incorporate the covariates $w$. Specifically, given $w$ and historical data $H_n$, a density estimation of the conditional distribution $\bbP[c \mid w]$ can be built using a kernel: $\bbP[\cdot \mid w] \approx \sum_{i \in [n]} k_{w_i}(w) \delta_{c_i}(\cdot)$ where $k_{w_i}(w)$ are convex combination weights which increase as the covariates $w$ become closer to $w_i$ (often obtained via a kernel), and $\delta_{c_i}(\cdot)$ is point mass at $c_i$. Then, a stochastic optimization problem with the estimated conditional distribution can be solved. This approach was studied by \citet{HannahPowellBlei2010,HanasusantoKuhn2013,BertsimasKallus2014,BanRudin2019,BertsimasVP2017,HoHanasusanto2019} who all gave various performance guarantees. 
{\cbl However, density estimation-based methods are known to require much more data than parametric prediction-based methods. As a result, when a reasonable parametric prediction model is available, it is advantageous to exploit it. Hence, density estimation methods are not the focus of this paper.}


\section{Risk Minimization and Consistency for Prediction and Optimization}\label{sec:jpo-risk}

{\cbl Given an vector $d$, recall from Section \ref{sec:jpo-intro} that we assess the quality of using $d$ in place of a true cost vector $c$ in \eqref{eqn:jpo-opt-problem} via the optimality gap of the solution obtained with $d$ on the true objective vector $c$, which we define to be the \emph{true loss function}
\begin{equation}\label{eqn:jpo-true-loss}
L(d,c) := f(x^*(d)) + c^\top x^*(d) - \min_{x \in X} \left\{ f(x) + c^\top x \right\},
\end{equation}
where $x^*(d)\in\argmin_{x \in X} \left\{ f(x) + d^\top x \right\}$ is as described in Notation subsection.} 
Note that given any $c \in \bbR^m$, $L(d,c) \geq 0$ for all $d \in \bbR^m$, and $L(c,c)=0$.
\begin{remark}
{\cbl By definition the true loss function} $L$ depends on the function $x^*$, i.e., the algorithm that we use to solve $\min_{x \in X} \left\{ f(x) + d^\top x \right\}$ for different $d \in \bbR^m$. We will take $x^*$ to be fixed throughout the paper. Note, however, that the specific choice of $x^*$ only affects our results up to measurability concerns; we show in Lemma \ref{lemma:jpo-any-x^*-measurable} that any $x^*$ is Lebesgue measurable, so we can safely fix $x^*$ without changing the results as long as our distribution $\bbP$ is Lebesgue measurable. In practice, any distribution we encounter will be Lebesgue measurable; we explicitly impose this in Assumption \ref{ass:jpo-lebesgue-measurability}. Henceforth, when measurability of functions is discussed, we will understand this to be in the sense of Lebesgue.
\epr
\end{remark}

{\cbl
Our setting of interest is prediction in the context of solving the optimization problem \eqref{eqn:jpo-opt-problem}. Specifically, instead of a solitary (random) cost vector $c$, we are interested in random pairs $(w,c) \in W \times \bbR^m$ drawn from a distribution $\bbP$. We are then interested in learning a prediction function $g:W \to \bbR^m$ which predicts $c$ with $g(w)$. Since $(w,c) \sim \bbP)$ is random, we assess the performance of the prediction function $g$ in terms of the expected true loss,  which we call the \emph{true risk}
\begin{equation}\label{eqn:jpo-true-risk}
R(g,\bbP) := \bbE[L(g(w),c)].
\end{equation}
The best possible true risk we can achieve is
\begin{equation}\label{eqn:jpo-true-risk-min}
R(\bbP):=\inf_g \left\{ R(g,\bbP) : g \text{ measurable} \right\}.
\end{equation}
}

A na\"{i}ve attempt to solve \eqref{eqn:jpo-true-risk-min} is to directly minimize $R(g,\bbP)$. However, as we show in {\cbl Lemma \ref{lemma:jpo-L-nonconvex} below,} $L(d,c)$ is not convex in $d$. Thus, in general it is not %
{\cbl expected} to obtain a polynomial-time approach to minimize the true risk $R(g,\bbP)$.
{\cbl
\begin{lemma}\label{lemma:jpo-L-nonconvex}
Suppose that $f(x) = 0$, $X$ is such that it has at least two extreme points and  $c\neq0$ is such that $\min_{x\in X}c^\top x \neq \max_{x\in X} c^\top x$. Then, the loss function $L(d,c)$ is not convex in $d$.
\end{lemma}
}

{\cbl Given Lemma~\ref{lemma:jpo-L-nonconvex}, in order to examine practical polynomial-time solution methodologies, we now describe an alternative approach based on \emph{surrogate loss functions}.} We first state a basic fact about \eqref{eqn:jpo-true-risk-min}, {\cbl which is analogous to \cite[Proposition 5]{ElmachtoubGrigas2017} adapted to our setup.}
\begin{lemma}\label{lemma:jpo-true-risk-minimizers}
The function $g^*(w) := \bbE[c \mid w]$ minimizes \eqref{eqn:jpo-true-risk-min}. Furthermore,
\[ R(\bbP) = \bbE\left[ \min_{d' \in \bbR^m} \bbE[L(d',c) \mid w] \right]. \]
\end{lemma}
{\ctl Lemma \ref{lemma:jpo-true-risk-minimizers} shows that the conditional expectation is a minimizer of \eqref{eqn:jpo-true-risk-min}. {\cbl There are several regression techniques which aim to recover the conditional expectation. These have a similar structure, which we now describe. First, we specify some loss function $\ell(d,c)$ measuring the discrepancy between vectors $d$ and $c$. As before, we are in the setting where we have random pairs $(w,c) \sim \bbP$ and we wish to do prediction via a function $g(w) \approx c$. We then define the \emph{surrogate risk} as:
\begin{equation}\label{eqn:jpo-surrogate-risk}
R_\ell(g,\bbP) := \bbE\left[ \ell(g(w),c) \right]
\end{equation}
as well as the best possible surrogate risk 
\begin{equation}\label{eqn:jpo-surrogate-risk-min}
R_\ell(\bbP) := \inf_g \left\{ R_\ell(g,\bbP) : g \text{ measurable} \right\}.
\end{equation}
Now, instead of seeking to minimize the true risk $R(g,\bbP)$, we seek to minimize the surrogate risk $R_\ell(g,\bbP)$. Indeed, regression methods often use tractable convex losses $\ell$, hence minimizing the surrogate risk is much more tractable than the true risk. Here, we use the term `surrogate' since, in a sense, the loss function $\ell$ can be thought of as a \emph{surrogate loss} for $L$, i.e., in order to maintain computational tractability, we replace the difficult loss $L$ with a more computationally friendly surrogate $\ell$.
}
The use of surrogate loss functions to ensure algorithmic tractability is very common in machine learning. For example, convex surrogates such as hinge loss are used instead of the non-convex true $0$-$1$ loss in classification problems.

It is not yet well-understood how minimizing the surrogate risk \eqref{eqn:jpo-surrogate-risk} can impact the true \eqref{eqn:jpo-true-risk}. A good surrogate loss function $\ell$ should mimic the natural properties of the true loss function $L$, i.e., $\ell(c,c)=0$, $\ell(d,c)\geq0$ for any $d,c$. However, the most important feature of a surrogate loss function is how its risk bound {\cbl relates} to the true risk \eqref{eqn:jpo-true-risk}. {\cbl More precisely, if one were to obtain a prediction function $\hat{g}$ with low excess surrogate risk $R_{\ell}(\hat{g},\bbP) - R_{\ell}(\bbP)$, will it be the case that $\hat{g}$ also has low excess true risk $R(\hat{g},\bbP) - R(\bbP)$?}

Consequently, in this paper we will explore this relationship and identify important properties of surrogate loss functions that enable us to derive guarantees on the true risk.
We would like to identify essential properties of surrogate loss functions $\ell(g(w),c)$ such that they can accurately, in some sense, assess the quality of using $g(w)$ in place of $c$ for the true risk \eqref{eqn:jpo-true-risk} related to \eqref{eqn:jpo-opt-problem}, while remaining computationally tractable to optimize (e.g., being convex in $g(w)$).

{\cbl While the concepts we explore are related to the more traditional notion of \emph{statistical consistency}, they are of a slightly different nature, which we elaborate on now. In practice, the distribution $\bbP$ is not given explicitly, but instead we only have access to historical data $H_n = \{(w_i,c_i) : i \in [n]\}$. To obtain a predictor $g:W \to \bbR^m$, we optimize the \emph{empirical} surrogate risk
\[\hat{R}_\ell(g,H_n) := \frac{1}{n} \sum_{i=1}^n \ell(g(w_i),c_i).\]
Statistical learning theory has rich literature on relating $\hat{R}_\ell$ to $R_\ell$; see, e.g., \citet{BousquetEtAl2004}. In particular, it has several results on the following notion of consistency.
\begin{definition}\label{def:jpo-statistical-consistency}
Given a (deterministically expanding) sequence of classes of predictors $\left\{\cG_n\right\}_{n \in \bbN}$, let $\hat{g}_n := \argmin_{g \in \CG_n} \hat{R}_\ell(g,H_n)$. We say that the (random) sequence of predictors $\{\hat{g}_n\}_{n \in \bbN}$ is \emph{statistically consistent with respect to loss $\ell$} if
\[ R_\ell(\hat{g}_n,\bbP) \to R_\ell(\bbP) \text{ in probability}. \]
\end{definition}
(Convergence in probability is used due to the randomness in $H_n$, which translates to randomness of $\hat{g}_n$.) This states that, for large $n$, we can get high-probability bounds on the excess surrogate risk $R_\ell(\hat{g}_n,\bbP) - R_\ell(\bbP)$ of a predictor $\hat{g}_n$. Whenever $\{\cG_n\}_{n \in \bbN}$ and $\ell$ are mildly regular, the consistency of the predictors $\hat{g}_n$ holds in a wide variety of settings.

However, since we will use $\hat{g}_n$ for optimization, we are actually interested in the excess true risk $R(\hat{g}_n,\bbP) - R(\bbP)$ which depends on \eqref{eqn:jpo-opt-problem} explicitly. Thus, in this paper, we give relationships between the excess surrogate risk and the true risk. More precisely, we will explore conditions on the surrogate loss function $\ell$ that ensure the following property holds:
\begin{definition}\label{def:jpo-L-consistency}
Given a class of distributions $\cP$, we say that $\ell$ is \emph{$(\cP,L)$-consistent} if, for all $\bbP \in \cP$, whenever we have a sequence of predictors $\{g_n\}_{n \in \bbN}$ such that $R_\ell(g_n, \bbP) \to R_\ell(\bbP)$, we will also imply $R(g_n,\bbP) \to R(\bbP)$.
\end{definition}
If $\ell$ satisfies Definition \ref{def:jpo-L-consistency} then this means that \emph{any} sequence that is statistically consistent (in the sense of Definition \ref{def:jpo-statistical-consistency}) with respect to the surrogate loss $\ell$ is also statistically consistent with respect to the true loss $L$.

\begin{remark}\label{rem:consistency-results}
Note that Definition \ref{def:jpo-L-consistency} does not depend on the classes of predictors $\{\CG_n\}_{n \in \bbN}$ or the sequence of predictors $\{\hat{g}_n\}_{n \in \bbN}$ obtained by minimizing $\hat{R}_\ell(g,H_n)$, even though these are important to relate the empirical surrogate risk $\hat{R}_\ell$ to the population surrogate risk $R_\ell$, as well as for computational considerations of optimizing the surrogate risk. Because of this, our results are also naturally independent of the choice of $\{\CG_n\}_{n \in \bbN}$. This is important for the application of our theory: by keeping the $\{\CG_n\}_{n \in \bbN}$ unspecified, our results are applicable to all settings.
\epr
\end{remark}
}
}


\section{Risk Minimization via Fisher Consistent Surrogate Loss Functions}\label{sec:jpo-FisherConsistent-surrogate-risk}

As discussed in Section \ref{sec:jpo-risk}, we are interested in properties of the surrogate loss $\ell$ which ensures consistency in the sense of Definition \ref{def:jpo-L-consistency} holds. In order to understand the kind of results that we are after, let us explore the negation of this. In this case, we have $R_\ell(g_n,\bbP) - R_\ell(\bbP) \to 0$ but, for some $\epsilon > 0$, $R(g_n,\bbP) - R(\bbP) > \epsilon$ for infinitely many $n$. In other words, there exists $\epsilon > 0$ such that for all $\delta > 0$, there exists $g_n$ such that $R_\ell(g_n,\bbP) - R_\ell(\bbP) \leq \delta$ but $R(g_n,\bbP) - R(\bbP) > \epsilon$. To prevent this bad outcome, we want to guarantee the following relationship between the risks:
\begin{align}
&\text{for all $\epsilon > 0$, there exists $\delta > 0$ such that:}\label{eqn:jpo-risk-relationship}\\
&\quad \text{if $g:W \to \bbR^m$ satisfies } R_\ell(g,\bbP) - R_\ell(\bbP) \leq \delta, \text{ then } R(g,\bbP) - R(\bbP) \leq \epsilon.\notag
\end{align}

We will show that \eqref{eqn:jpo-risk-relationship} can be guaranteed by checking a simpler condition on the losses $\ell$ and $L$ called \emph{calibration}. This was introduced by \citet{BartlettJordanMcAuliffe2006} for binary classification and extended by \citet{Steinwart2007} for other machine learning applications. We extend this concept to the context of prediction and optimization.
\begin{definition}\label{def:jpo-calibration}
A surrogate loss function $\ell$ for $L$ is \emph{calibrated} with respect to a distribution $\bbP$, or \emph{$\bbP$-calibrated}, if, for all $w \in W$ and $\epsilon > 0$, there exists $\delta > 0$ (which may depend on $w$) such that
\[ \text{if $d \in \bbR^m$ satisfies } \bbE[\ell(d,c) \mid w] - \min_{{\cbl d'} \in \bbR^m} \bbE[\ell({\cbl d'},c) \mid w] ~{\cbl <}~ \delta, \text{ then } \bbE[L(d,c) \mid w] - \min_{{\cbl d'} \in \bbR^m} \bbE[L({\cbl d'},c) \mid w] ~{\cbl <}~ \epsilon. \]
\end{definition}

Observe that Definition \ref{def:jpo-calibration} is very similar to \eqref{eqn:jpo-risk-relationship}, except that predictors $g$ (i.e., functions mapping onto vectors) are replaced with vectors $d \in \bbR^m$, {\cbl and that conditional expectations given $w$ are taken. This makes Definition \ref{def:jpo-calibration} verifiable, i.e., given a class of probability distributions $\bbP$ and a surrogate loss $\ell$, we can check whether Definition \ref{def:jpo-calibration} holds or not. Of course, we then need to check that Definition \ref{def:jpo-calibration} is sufficient to obtain risk bounds. \citet[Theorem 2.8]{Steinwart2007} provides a result to obtain such risk bounds, and we apply it to obtain Theorem \ref{thm:jpo-calibration} below.
More precisely, we verify that necessary measurability and boundedness conditions on certain conditional risk quantities are met in order to apply the proof technique of \citet{Steinwart2007} in our prediction and optimization context.}
We use the following technical assumption:
\begin{assumption}\label{ass:jpo-finiteness}
Let the probability distribution $\bbP$ and the surrogate loss function $\ell$ be given. For any fixed $c \in C$, the surrogate loss function $\ell(d,c)$ is convex in $d \in \bbR^m$. For any $w \in W$ and $d \in \bbR^m$, the set $\argmin_{d' \in \bbR^m} \bbE[\ell(d',c) \mid w]$ is non-empty and bounded, and $\bbE[\ell(d,c) \mid w] < \infty$. Furthermore, $c$ is an integrable random vector (that is, each component is integrable) so that $\bbE[\|c\|_1] < \infty$.
\end{assumption}

\begin{theorem}\label{thm:jpo-calibration}
Suppose that $\ell$ is $\bbP$-calibrated, and that Assumption \ref{ass:jpo-finiteness} holds. Then for all $\epsilon > 0$, there exists a $\delta > 0$ such that
\[ R_\ell(g,\bbP) \leq R_\ell(\bbP) + \delta \implies R(g,\bbP) \leq R(\bbP) + \epsilon. \]
\end{theorem}
We give the proof in Section \ref{sec:jpo-proof-calibration}.

In general, checking that a given surrogate loss $\ell$ is $\bbP$-calibrated {\cbl may not be straightforward}. A much simpler condition to check is Fisher consistency, %
{\cbl stated} in Definition \ref{def:jpo-FisherConsistent-loss} below. %
{\cbl Note that Fisher consistency}  relates to the \emph{minimizers} of the loss functions, instead of approximate minimizers as in Definition \ref{def:jpo-calibration}. In this section, we show that Fisher consistency is equivalent to calibration, thus allowing us to check the simpler condition to verify Theorem \ref{thm:jpo-calibration}. We also discuss some different loss functions and their Fisher consistency properties.
\begin{definition}\label{def:jpo-FisherConsistent-loss}
A surrogate loss function $\ell$ is \emph{Fisher consistent} with respect to a distribution $\bbP$, or \emph{$\bbP$-Fisher consistent}, if for all $w$,
\[ \argmin_d \bbE\left[ \ell(d,c) \mid w \right] \subseteq \argmin_d \bbE\left[ L(d,c) \mid w \right]. \]
\end{definition}
{\cbl
\begin{remark}
Since the objective for our optimization problem is of the form $f(x) + c^\top x$, we proved in Lemma \ref{lemma:jpo-true-risk-minimizers} that $\bbE[c \mid w] \in \argmin_d \bbE\left[ L(d,c) \mid w \right]$. Thus, one way to check that $\ell$ is Fisher consistent is to verify that $\argmin_d \bbE\left[ \ell(d,c) \mid w \right] = \left\{ \bbE[c \mid w] \right\}$ (and this is the approach taken in some of the examples below). However, we opt not to simply take $\argmin_d \bbE\left[ \ell(d,c) \mid w \right] = \left\{ \bbE[c \mid w] \right\}$ as the definition of Fisher consistency because we recognize, particularly for non-smooth optimization objectives, that there can be other vectors besides $\bbE[c \mid w]$ that minimize $\bbE\left[ L(d,c) \mid w \right]$. Furthermore, the current form of Definition \ref{def:jpo-FisherConsistent-loss} will also allow us to encompass settings when the objective is of a more general form than $f(x) + c^\top x$ (although this is not the focus of the current paper).
\epr
\end{remark}
}

{\cbl In} the following theorem {\cbl we show} that Fisher consistency is equivalent to calibration. Of course, the fact that calibration implies Fisher consistency is straightforward; the {\cbl main} challenge is to show the other direction.
\begin{theorem}\label{thm:jpo-FisherConsistent-calibration}
Given a distribution $\bbP$, let $\ell(d,c)$ be {\cbl a} loss function that satisfies Assumption \ref{ass:jpo-finiteness}. Then $\ell$ is $\bbP$-calibrated if and only if $\ell$ is $\bbP$-Fisher consistent.
\end{theorem}
The key tool that we exploit {\cbl in proving Theorem~\ref{thm:jpo-FisherConsistent-calibration}} is upper semi-continuity of the multivalued argmin mapping $X^*(\cdot)$ (see \ref{sec:jpo-regularity}). Informally, this states that if we are given $X^*(d)$ for some vector $d$, and we are interested in vectors $d'$ for which $X^*(d')$ does not move `too far away' from $X^*(d)$, then we can guarantee that when $d'$ is sufficiently close to $d$, this will indeed be the case. In particular, in the context of proving Theorem \ref{thm:jpo-FisherConsistent-calibration}, we use this to show that when $\bbE[L(d,c) \mid w]$ is large, then vectors close by to $d$ will also have large true expected loss. The full proof of Theorem \ref{thm:jpo-FisherConsistent-calibration} is in \ref{sec:jpo-proofs-FisherConsistent}.

Armed with Theorem \ref{thm:jpo-FisherConsistent-calibration}, we have the following corollaries, which are straightforward consequences of {\cbl our} %
results {\cbl discussed so far}.
\begin{corollary}\label{cor:jpo-FisherConsistent}
Suppose that $\ell$ is $\bbP$-Fisher consistent, and that Assumption \ref{ass:jpo-finiteness} holds. Then for all $\epsilon > 0$, there exists a $\delta > 0$ such that
\[ R_\ell(g,\bbP) \leq R_\ell(\bbP) + \delta \implies R(g,\bbP) \leq R(\bbP) + \epsilon. \]
\end{corollary}

\begin{corollary}\label{cor:jpo-FisherConsistent-conv}
Suppose that $\ell$ is $\bbP$-Fisher consistent, and that Assumption \ref{ass:jpo-finiteness} holds. If we have a sequence of functions $g_n$ such that $R_\ell(g_n,\bbP) \to R_\ell(\bbP)$. Then $R(g_n,\bbP) \to R(\bbP)$.
\end{corollary}

We now examine several different loss functions and their Fisher consistency properties. Before doing so, let us summarize the properties on $\ell$ and $\bbP$ in order to get risk guarantees of the form \eqref{eqn:jpo-risk-relationship} {\cbl through Theorem~\ref{thm:jpo-FisherConsistent-calibration}}. These are:
\begin{enumerate}
\item\label{item:jpo-convex} the surrogate loss $\ell(\cdot,c)$ is convex for any fixed $c \in C$.
\item\label{item:jpo-finite} for any $w \in W, d \in \bbR^m$, the expected loss $\bbE[\ell(d,c) \mid w]$ is finite.
\item\label{item:jpo-bounded-minimizers} for any $w \in W$, the set of minimizers $\argmin_{d' \in \bbR^m} \bbE[\ell(d',c) \mid w]$ is non-empty and bounded.
\item\label{item:jpo-FisherConsistent} the surrogate loss $\ell$ is $\bbP$-Fisher consistent according to Definition \ref{def:jpo-FisherConsistent-loss}.
\end{enumerate}

We first examine {\cbl the squared loss function, namely $\ell_{\LS}(d,c) = \|d - c\|_2^2$,} that is Fisher consistent for any class of distributions. {\cbl (We use the `$\LS$' subscript as shorthand for `least squares'.)}
\begin{example}\label{ex:jpo-squaredLoss}
Consider the squared loss $\ell_{\cbl \LS}(d,c) = \|d - c\|_2^2$. Then $\ell_{\cbl \LS}$ is $\bbP$-Fisher consistent for any distribution $\bbP$ over $W \times C$. Note that 
\[ \bbE[\ell_{\cbl \LS}(d,c) \mid w] = \bbE\left[ \|d - c\|_2^2 {\cbl \mid w} \right] = \|d - \bbE[c \mid w]\|_2^2 + \bbE[\|c\|_2^2 \mid w] - \|\bbE[c \mid w]\|_2^2. \]
Thus, the unique minimizer of $\bbE[\ell_{\cbl \LS}(d,c) \mid w]$ is $d^* = \bbE[c \mid w]$. Since we know this is also a minimizer of $\bbE[L(d,c) \mid w]$, this gives us $\bbP$-Fisher consistency of the squared loss; {\cbl verifying Property~\ref{item:jpo-FisherConsistent}.}

Also note that {\cbl Properties} \ref{item:jpo-convex} and \ref{item:jpo-bounded-minimizers} are clearly satisfied. Property \ref{item:jpo-finite} will be satisfied if the conditional distribution $\bbP[\cdot \mid w]$ is square integrable for every $w \in W$.
\epr
\end{example}

A common loss function used in regression to safeguard against outliers is the absolute deviation loss, {\cbl namely $\ell_{\cbl \AD}(d,c) := \|d - c\|_1$. We next examine this loss function.}
\begin{example}\label{ex:jpo-abs-dev-loss}
Consider the absolute deviation loss $\ell_{\cbl \AD}(d,c) = \|d - c\|_1$. %
We claim that $\ell_{\cbl \AD}$ is $\bbP$-Fisher consistent as long as, for every $w$, $\bbP[\cdot \mid w]$ is centrally symmetric about some vector $d_w$. A distribution $\bbP$ is centrally symmetric about $d$ if, for a random variable $c \sim \bbP$, $c-d$ has the same distribution as $d-c$. Note that $\argmin_{d' \in \bbR^m} \bbE[\|d'-c\|_1 \mid w]$ recovers the vector of coordinate-wise medians, which for a centrally symmetric distribution will be the point of symmetry $d_w$, i.e., $d_w$ minimizes $\bbE[\|d-c\|_1 \mid w]$. Furthermore, we have $\bbE[c \mid w] = d_w$ also. 
Therefore, $d_w$ minimizes $\bbE[L(d,c) \mid w]$.
\epr
\end{example}

We now discuss the $\SPOp$ loss function proposed in \citet{ElmachtoubGrigas2017}, which aims to incorporate knowledge of the domain $X$ into the loss, in the hopes of achieving low true risk $R$, which is based on the optimization problem.
\begin{example}\label{ex:jpo-spo-loss}
In the setting when $f(x) = 0$ for all $x \in X$, \citet[Definition 3]{ElmachtoubGrigas2017} defined the following loss function:
\begin{equation}\label{eqn:jpo-spo-loss}
\ell_{\SPOp}(d,c) := (2d-c)^\top x^*(c) - \min_{x \in X} (2d-c)^\top x = L(c,2d-c).
\end{equation}
\citet[Proposition 6]{ElmachtoubGrigas2017} shows that $\ell_{\SPOp}$ is Fisher consistent as long as $\bbP[c \mid w]$ is centrally symmetric and continuous. We remark also that \citet{ElmachtoubGrigas2017} achieve good numerical results, particularly when the hypothesis class is misspecified versus the true distribution.
\epr
\end{example}

We now highlight some positive and negative aspects of the loss function of \citet{ElmachtoubGrigas2017}. {\cbl We start with an example below to review} an important observation made in \citet[Proposition 1]{ElmachtoubGrigas2017} {\cbl that in the case of binary classification,} by carefully choosing the set $C$ and domain $X$, the true loss $L$ from \eqref{eqn:jpo-true-loss} becomes the $0$-$1$ loss. {\cbl In addition, their surrogate} loss {\cbl $\ell_{\SPOp}$} \eqref{eqn:jpo-spo-loss} also has a familiar interpretation {\cbl and admits Fisher consistency} in this setting. 
\begin{example}\label{ex:jpo-binary-class}
Let $m=1$, $C=\{-1,1\}$, $X=[-1/2,1/2]$ and $f(x) = 0$ for all $x \in X$. Then, $x^*(d) = -\sign(d)/2$, and $\min_{x \in X} c^\top x = -1/2$ for any $c \in C$, so
\[ L(d,c) = \frac{c \sign(d) - 1}{2} = \begin{cases} 0, & c = \sign(d)\\ 1, & c \neq \sign(d). \end{cases} \]
{\cbl That is,} {\ctl the $0$-$1$ loss for classification is exactly equivalent to the true loss function $L$.}
\citet[Proposition 4]{ElmachtoubGrigas2017} shows that the loss from \eqref{eqn:jpo-spo-loss} reduces to the hinge loss in this case: since $x^*(c) = -c/2$ for $c \in C$ and $\min_{x \in X} d^\top x = -|d|/2$,
\begin{align*}
\ell_{\SPO+}(d,c) &= \frac{|2d-c| - (2d-c)c}{2} = \frac{|1-2dc| + 1-2dc}{2} = \max\{0,1-2dc\} . 
\end{align*}
{\cbl Moreover,} \citet[Theorem 3.1]{Lin2004} states that the hinge loss, and thus $\ell_{\cbl \SPO+}$, is Fisher consistent for any distribution over $C = \{-1,1\}$ except the uniform one.
\epr
\end{example}

{\cbl In contrast to this, we next demonstrate with the following two general examples that} %
the loss function {\cbl $\ell_{\SPO+}$} of \citet{ElmachtoubGrigas2017} is not Fisher consistent %
{\cbl in} some very natural settings. %
\begin{example}\label{ex:jpo-spo-one-dim}
Consider the setting where $m=1$, $X = [-1/2,1/2]$ and $f(x)=0$ for all $x \in X$, but $C$ is an arbitrary subset of $\bbR$. Then $x^*(c) = -\sign(c)/2$, $\min_{x \in X} d^\top x = -|d|/2$, hence the loss function from \eqref{eqn:jpo-spo-loss} becomes
\[ \ell_{\SPOp}(d,c) = \frac{|2d-c| - (2d-c)\sign(c)}{2} = \frac{|2d-c| - 2d\sign(c) + |c|}{2}. \]
Let $\bbP$ be a distribution over $W \times C$. For any $w \in W$, note that the minimizers of $\bbE[L(d,c) \mid w]$ are $D_w^* = \{d \in \bbR : \sign(d) = \sign(\bbE[c \mid w])\}$. Thus, checking $\bbP$-Fisher consistency requires showing that $\argmin_{d' \in \bbR} \bbE[\ell_{\SPO+}(d',c) \mid w] \subseteq D_w^*$ for every $w \in W$, i.e., we need to show that the minimizers have the same sign as the mean $\bbE[c \mid w]$. However, we %
show (in \ref{sec:jpo-proofs-FisherConsistent}) that the minimizer of the loss function {\cbl $\ell_{\SPOp}(d,c)$} %
has the same sign as the median. {\cbl Therefore,} %
for distributions where the mean and median have different signs, this loss function is not Fisher consistent.
\epr
\end{example}

\begin{example}\label{ex:jpo-multiclass}
In Example \ref{ex:jpo-binary-class}, {\cbl we examined binary classification and} showed that for appropriately chosen $X$, $f$ and $C$, $L$ specializes to the $0$-$1$ loss and $\ell_{\SPOp}$ specializes to the hinge loss. Thus, $\ell_{\SPOp}$ defined in \eqref{eqn:jpo-spo-loss} can be seen as a generalization of the hinge loss for optimization problems. We {\cbl next} show that %
the multiclass classification loss {\cbl admits a similar representation}, i.e., {\cbl by} choosing $X$ and $C$ appropriately we can make $L$ represent the $0$-$1$ loss for multiclass classification. However, {\cbl we also establish that} the generalization of hinge loss given by \eqref{eqn:jpo-spo-loss} to this setting is not Fisher consistent.

Suppose we have pairs $(w,c)$, where $w$ are features, and $c \in C'$ is a label from one of $m \in \bbN$ different classes, i.e., $C' = [m]$. We want a predictor $g':W \to C'$ which classifies $w$ according to $g'(w)$. If we classify $w$ incorrectly (i.e., $g'(w)$ is in a different class to $c$) we suffer {\cbl a} loss {\cbl of} $1$; otherwise, our loss is $0$. We can capture this in our optimization framework as follows. 

Consider $C = \{c_j := \bm{1}_m - e_j : j \in [m]\} \subset \bbR^m$, $X = \Conv\left\{ e_j : j \in [m] \right\} \subset \bbR^m$ and $f(x) = 0$ for all $x \in X$. Then $\min_{x \in X} d^\top x = \min_{j' \in [m]} d_{j'}$, $\min_{x \in X} c_j^\top x = 0$ and $x^*(d) = e_j$ for $j \in \argmin_{j' \in [m]} d_{j'}$, so for any $j \in [m]$ and vector $d$ with unique minimum entry
\[ L(d,c_j) = \begin{cases}
0, & \argmin_{j' \in [m]} d_{j'} = j\\
1, & \argmin_{j' \in [m]} d_{j'} \neq j.
\end{cases} \]
In other words, if we have a function $g:W \to \bbR^m$, we can use it to build a classifier $g':W \to C'$ by classifying $w$ according to the minimum entry of $g(w) \in \bbR^m$. Then $L$ is exactly the $0$-$1$ loss for this classifier. Suppose that we have a distribution $\bbP[c=c_j] = p_j > 0$, $\sum_{j \in [m]} p_j = 1$. Then, letting $j^*(d) = \argmin_{j' \in [m]} d_{j'}$,
\[ \bbE[L(d,c)] = 1 - p_{j^*(d)}, \]
so the vectors $d$ which minimize $\bbE[L(d,c)]$ must satisfy $j^*(d) \in \argmax_{j' \in [m]} p_{j'}$.

The loss \eqref{eqn:jpo-spo-loss} becomes
\[ \ell_{\SPOp}(d,c_j) = (2d-c_j)^\top e_j - \min_{j' \in [m]} \left\{ 2d_{j'} - c_{j'} \right\} = 2d_j - \min_{j' \in [m]} \left\{ 2d_{j'} - \bm{1}(j' \neq j) \right\}. \]
In \ref{sec:jpo-proofs-FisherConsistent}, we show that for distributions $\bbP$ with $\max_{j' \in [m]} p_{j'} < 1/2$, $\ell$ is not $\bbP$-Fisher consistent, since the set of minimizers of $\bbE[\ell_{\SPOp}(d,c)]$ are the vectors $d_\alpha = \alpha \bm{1}_m$, $\alpha \in \bbR$, which cannot in general pick out the maximum probability class $j \in [m]$, i.e., the highest $p_j$.
\epr
\end{example}


\section{Non-Asymptotic Risk Guarantees via Uniform Calibration}\label{sec:jpo-non-asymptotic}

Corollary \ref{cor:jpo-FisherConsistent-conv} is an asymptotic result, that is, %
{\cbl it asserts only} that minimizing the surrogate risk will minimize the true risk in the limit. This does not present much insight about the rate of convergence of these quantities, which is governed by the relationship between $\epsilon$ and $\delta$ in Corollary \ref{cor:jpo-FisherConsistent}. Moreover, the $\delta$ in Corollary \ref{cor:jpo-FisherConsistent} depends on the distribution $\bbP$. In general, this is undesirable, since often in statistical learning, we assume minimal knowledge of $\bbP$. Furthermore, when given $n$ data points $\{(w_i,c_i) : i \in [n]\}$ we can build a predictor $g_n$ with quantified guarantees on the excess surrogate risk $R_\ell(g_n,\bbP) - R_\ell(\bbP)$ via standard learning theoretic results. We would ideally like to translate these into quantified guarantees on the excess true risk $R(g_n,\bbP) - R(\bbP)$.

\citet{Steinwart2007} builds a theory for non-asymptotic relationships between true and surrogate risk for various types of learning problems, such as classification, regression, and density estimation, giving necessary and sufficient conditions for the existence of distribution-independent guarantees. In this section, building on the results from \citet{Steinwart2007}, we provide conditions for the existence of similar guarantees in the prediction and optimization context. Using these conditions, we identify a non-asymptotic distribution-independent guarantee between the risk of the surrogate squared loss function {\cbl $\ell_{\LS}$} and the true {\cbl optimality} gap risk. We then provide risk guarantees for a class of symmetric loss functions by appealing to existing results on their risk relationships to the squared loss {\cbl $\ell_{\LS}$}. Finally, we study the a special case of the  {\cbl $\ell_{\SPOp}$} loss function \eqref{eqn:jpo-spo-loss} of \citet{ElmachtoubGrigas2017}, and provide positive and negative results on %
{\cbl its} risk guarantees.

\subsection{Outline of %
Key Idea}\label{sec:jpo-non-asymptotic-outline}

{\cbl In order to provide guarantees on the true risk implied by the surrogate risk,} in this section, our aim is to identify an increasing function $\eta:[0,\infty) \to [0,\infty)$ with $\eta(0) = 0$ such that for any distribution $\bbP$, we have
\[ \eta\left( R(g,\bbP) - R(\bbP) \right) \leq R_\ell(g,\bbP) - R_\ell(\bbP). \]
Thus, any bound on the excess surrogate risk $R_\ell(g,\bbP) - R_\ell(\bbP)$ translates to a bound on the excess true risk $R(g,\bbP) - R(\bbP)$. Let us explore how we would derive such bounds. First, suppose that {\cbl $\eta$ and $\ell$ are chosen so that $\eta$ is convex and that for any $w \in W$ and $d \in \bbR^m$, we have
\begin{equation}\label{eqn:jpo-eta-relationship}
		\eta\left( \bbE\left[ L(d,c) \mid w \right] - \min_{d' \in \bbR^m} \bbE\left[ L(d',c) \mid w \right] \right) \leq \bbE\left[ \ell(d,c) \mid w \right] - \min_{d' \in \bbR^m} \bbE\left[ \ell(d',c) \mid w \right].
\end{equation}
Then, we have
\begin{align*}
\eta\left( R(g,\bbP) - R(\bbP) \right) &= \eta\left( \bbE\left[ \bbE\left[ L(g(w),c) \mid w \right] - \min_{d' \in \bbR^m} \bbE\left[ L(d',c) \mid w \right] \right] \right)\\
&\leq \bbE\left[ \eta\left( \bbE\left[ L(g(w),c) \mid w \right] - \min_{d' \in \bbR^m} \bbE\left[ L(d',c) \mid w \right] \right) \right]\\
&\leq \bbE\left[ \bbE\left[ \ell(g(w),c) \mid w \right] - \min_{d' \in \bbR^m} \bbE\left[ \ell(d',c) \mid w \right] \right]\\
&= R_\ell(g,\bbP) - R_\ell(\bbP),
\end{align*}
where the first inequality follows from Jensen's inequality, and the second follows from \eqref{eqn:jpo-eta-relationship}.}

As a first attempt to choose such $\eta$ and $\ell$, we define
\begin{equation}\label{eqn:jpo-calibration-function-0}
\delta_\ell(\epsilon,w;\bbP) := \inf_{d \in \bbR^m} \left\{ \bbE[\ell(d,c) \mid w] - \min_{d' \in \bbR^m} \bbE[\ell(d',c) \mid w] :~ \bbE[L(d,c) \mid w] - \min_{d' \in \bbR^m} \bbE[L(d',c) \mid w] > \epsilon \right\}.
\end{equation}
{\cbl
\begin{remark}\label{rem:delta-intuition}
Note that $\delta_\ell(\epsilon,w;\bbP)$ is simply giving an explicit representation of the $\delta$ that appears in Definition \ref{def:jpo-calibration} as a function of $\epsilon$ and $w$. In particular, $\delta_{\ell}(\epsilon,w;\bbP) > 0$ for $\epsilon > 0$ whenever $\ell$ is $\bbP$-calibrated.

To see this, suppose that $\ell$ is $\bbP$-calibrated. Fix some $w \in W$ and $\epsilon > 0$. Then (by the contrapositive statement of the implication in Definition \ref{def:jpo-calibration}) there exists $\delta > 0$ such that whenever $\bbE[L(d,c) \mid w] - \min_{d' \in \bbR^m} \bbE[L('d,c) \mid w] > \epsilon$ we have $\bbE[\ell(d,c) \mid w] - \min_{d' \in \bbR^m} \bbE[\ell('d,c) \mid w] > \delta$. Taking the infimum of $\bbE[\ell(d,c) \mid w] - \min_{d' \in \bbR^m} \bbE[\ell('d,c) \mid w]$ over $d$ such that $\bbE[L(d,c) \mid w] - \min_{d' \in \bbR^m} \bbE[L('d,c) \mid w] > \epsilon$ gives exactly $\delta_{\ell}(\epsilon,w;\bbP)$ as defined in \eqref{eqn:jpo-calibration-function-0}, and we know that $\bbE[\ell(d,c) \mid w] - \min_{d' \in \bbR^m} \bbE[\ell('d,c) \mid w] > \delta$ for such $d$, hence $\delta_\ell(\epsilon,w;\bbP) \geq \delta > 0$.
\epr
\end{remark}
}
Fixing $w \in W$, consider $d \in \bbR^m$ such that $\bbE\left[ L(d,c) \mid w \right] - \min_{d' \in \bbR^m} \bbE\left[ L(d',c) \mid w \right] = \epsilon$. Then
\begin{align*}
\bbE\left[ \ell(d,c) \mid w \right] - \min_{d' \in \bbR^m} \bbE\left[ \ell(d',c) \mid w \right] &\geq \delta_\ell(\epsilon,w; \bbP)\\
&= \delta_\ell\left( \bbE\left[ L(d,c) \mid w \right] - \min_{d' \in \bbR^m} \bbE\left[ L(d',c) \mid w \right], w; \bbP \right).
\end{align*}
{\cbl This relation then inspires us to select $\eta = \delta_\ell$. But, such a choice of $\eta=\delta_\ell$ may not be feasible as we cannot ensure that $\delta_{\ell}$ is convex in general. Instead, we can} use $\eta = \delta_{\ell}^{**}$, where, given a function $h:\bbR \to \bbR \cup \{\infty\}$,
\begin{align*}
h^{**}(\epsilon) &= \sup_{h'} \left\{ h'(\epsilon) : h' \text{ convex function on $\bbR$}, h' \leq h \text{ pointwise} \right\}.
\end{align*}
Clearly, $h^{**}$ is convex since it is a supremum of convex functions, and it can be obtained via convex conjugacy (however, we will not need to appeal to this representation for our results).

Note that $\delta_{\ell}$ is only defined for $\epsilon > 0$, so we define $\delta_\ell({\cbl \epsilon},w;\bbP) = 0$ {\cbl when $\epsilon = 0$} and ${\cbl \delta_\ell({\cbl \epsilon},w;\bbP)} =+\infty$ when $\epsilon < 0$. Using $\eta = \delta_\ell^{**}$ guarantees {\cbl both} convexity {\cbl of $\eta$} and also that $\eta(\epsilon,w;\bbP) \leq \delta_\ell(\epsilon,w;\bbP)$, hence the desired inequality \eqref{eqn:jpo-eta-relationship} holds. Now, by the definition \eqref{eqn:jpo-calibration-function-0}, we have $\delta_\ell$ is non-decreasing in $\epsilon$ and positive for $\bbP$-calibrated $\ell$. However, $\ell$ could be such that $\delta_\ell(\epsilon,w;\bbP)$ does not increase once $\epsilon$ is sufficiently large, or only increases at a sublinear rate; in this case $\eta = \delta_\ell^{**}$ is going to be $0$ for $\epsilon \geq 0$, so the inequality \eqref{eqn:jpo-eta-relationship} will be useless. To prevent this, we make the assumption that $\bbE\left[ L(d,c) \mid w \right] - \min_{d' \in \bbR^m} \bbE\left[ L(d',c) \mid w \right] \leq B$ for all $w \in W$, $d \in \bbR^m$. We can then re-define $\delta_\ell(\epsilon,w;\bbP) = \infty$ for $\epsilon > B$, and take $\eta = \delta_{\ell}^{**}$. This ensures that $\eta(\epsilon) > 0$ for $\epsilon \in (0,B]$. To ensure that such a $B$ exists, we define the following quantities:
\begin{equation}\label{eq:jpo-uniform-diameters}
B_X := \max_{x,x' \in X} \|x-x'\|_2, \quad B_f := \max_{x,x' \in X} \left\{ f(x) - f(x') \right\}, \quad B_C := \max_{c \in C} \|c\|_2.
\end{equation}
Note that since $X$ is compact and $f$ is continuous on $X$, $B_X,B_f <\infty$.
\begin{assumption}\label{ass:jpo-C-compact}
The quantity $B_C < \infty$. (This means that $\bbE[c \mid w] \in \Conv(C)$ is uniformly bounded over $w \in W$.)
\end{assumption}
\begin{remark}\label{rem:jpo-boundedness}
Under Assumption \ref{ass:jpo-C-compact} and {\cbl using} the fact that $X$ is compact, we have
\begin{align*}
&\bbE\left[ L(d,c) \mid w \right] - \min_{d' \in \bbR^m} \bbE\left[ L(d',c) \mid w \right]\\
&= f(x^*(d)) - f(x^*(\bbE[c \mid w])) + \bbE[c \mid w]^\top \left( x^*(d) - x^*(\bbE[c \mid w]) \right)\\
&\leq f(x^*(d)) - f(x^*(\bbE[c \mid w])) + \left\| \bbE[c \mid w] \right\|_2 \left\| x^*(d) - x^*(\bbE[c \mid w]) \right\|_2\\
&\leq B_f + B_C B_X < \infty,
\end{align*}
{\cbl where the first inequality follows from Cauchy-Schwarz.}
\epr
\end{remark}

Another subtlety that we need to consider is that there needs to be {\cbl a \emph{single}} fixed $\eta$ for which \eqref{eqn:jpo-eta-relationship} holds for all $w \in W$. Because of this, the definition $\eta = \delta_\ell^{**}$ is not well-defined as $\delta_\ell$ in \eqref{eqn:jpo-calibration-function-0} depends on $w \in W$. To get around this, we need to strengthen the definition of calibration to be uniform across $w \in W$. In summary, the additions we need to make to the assumptions from Section \ref{sec:jpo-FisherConsistent-surrogate-risk} are Assumption \ref{ass:jpo-C-compact}, which ensures a uniform bound on the expected true loss, and a stronger definition of calibration, which we give next. Notice, however, that since {\cbl our} %
proof technique is different to that of Theorem \ref{thm:jpo-calibration}, we need {\cbl only} measurability of $\ell$, and not {\cbl necessarily its} convexity in $d$. In practice, however, convexity of $\ell$ in $d$ gives us implementable algorithms with performance guarantees.

\subsection{Risk Bounds via Uniform Calibration}

We consider the following strengthening of Definition \ref{def:jpo-calibration}.
\begin{definition}\label{def:jpo-uniform-calibration}
We say that a loss function $\ell$ is \emph{uniformly calibrated} with respect to a class of distributions $\CP$ on $W \times C$, or \emph{$\CP$-uniformly calibrated}, if, for all $\epsilon > 0$, there exists $\delta > 0$ such that for all $\bbP \in \CP$, $w \in W$ and $d \in \bbR^m$, we have
\begin{equation}\label{eqn:jpo-uniform-calibration}
\bbE[\ell(d,c) \mid w] - \inf_{d'} \bbE[\ell(d',c) \mid w] ~{\cbl <}~ \delta \implies \bbE[L(d,c)] - \inf_{d'} \bbE[L(d',c)] ~{\cbl <}~ \epsilon.
\end{equation}
\end{definition}
Note that Definition \ref{def:jpo-uniform-calibration} considers a class of distributions $\CP$ so that we can get distribution-independent guarantees. {\cbl This is due to practical considerations where knowledge of $\bbP$ may not be available explicitly, but rather} %
we may know that $\bbP$ belongs to some class $\CP$, so we may aim to get guarantees on the class $\CP$.

If $\ell$ is $\CP$-uniformly calibrated, then we define
\begin{equation}\label{eqn:jpo-uniform-calibration-function}
\delta_\ell(\epsilon;\CP) := \inf\limits_{\substack{d \in \bbR^m\\ w \in W\\ \bbP \in \CP}} \left\{ \bbE[\ell(d,c) \mid w] - \min\limits_{d' \in \bbR^m} \bbE[\ell(d',c) \mid w] :~ \bbE[L(d,c) \mid w] - \min\limits_{d' \in \bbR^m} \bbE[L(d',c) \mid w] ~{\cbl \geq}~ \epsilon \right\}.
\end{equation}
\begin{remark}\label{rem:jpo-uniform-calibration-function}
If $\ell$ is $\CP$-calibrated, then $\delta(\epsilon;\CP) > 0$ for all $\epsilon > 0$ {\cbl by taking the contrapositive of \eqref{eqn:jpo-uniform-calibration}}, and is non-decreasing in $\epsilon$. 
In addition, if Assumption \ref{ass:jpo-C-compact} holds, then $\delta_\ell(\epsilon;\CP) = \infty$ for $\epsilon > B_f + B_C B_X$ since the infimum is infeasible. Also, $\delta_\ell(\epsilon;\CP) = 0$ for $\epsilon < 0$. Furthermore, measurability of $\delta_\ell(\cdot;\CP)$ follows by a similar proof to Lemma \ref{lemma:jpo-delta-measurable}.
\epr
\end{remark}
Remark \ref{rem:jpo-uniform-calibration-function} shows that positivity of $\delta_\ell$ is necessary for $\CP$-uniform calibration. {\cbl We next establish} %
that it is also sufficient.
\begin{lemma}\label{lemma:jpo-uniform-calibration}
A surrogate loss function $\ell$ is $\CP$-uniformly calibrated {\cbl if and only if} $\delta_\ell(\epsilon;\CP) > 0$ for all $\epsilon > 0$. 
\end{lemma}

{\cbl We now have the tools to prove the risk guarantee for uniform calibration. This is presented as Theorem~\ref{thm:jpo-uniform-calibration} below, and we utilize a result of \citet[Theorem 2.13]{Steinwart2007} to prove it. Remark \ref{rem:jpo-boundedness} allows us to apply this result in the prediction and optimization context.
In this proof, it is crucial to ensure that the risk guarantee is non-trivial, i.e., verifying that $\delta_\ell^{**}$ is positive on its domain. We utilize Lemma \ref{lemma:jpo-uniform-calibration} for this purpose.}
\begin{theorem}\label{thm:jpo-uniform-calibration}
Suppose that $\ell$ is $\CP$-uniformly calibrated, and that Assumption \ref{ass:jpo-C-compact} holds. Define
\[\delta_\ell^{**}(\epsilon; \cP) := \sup_{h'} \left\{ h'(\epsilon) :~ h' \text{ convex function on $\bbR$}, ~h' \leq \delta_\ell(\cdot; \cP) \text{ pointwise on $(0,B_f + B_C B_X]$} \right\},\]
where $B_f, B_C, B_X$ are defined in \eqref{eq:jpo-uniform-diameters}. Then $\delta_\ell^{**}(\epsilon;\CP)$ is positive for $\epsilon \in (0,B_f + B_C B_X]$, and for any $\bbP \in \CP$, $g:W \to \bbR^m$,
\[ \delta_\ell^{**}\left( R(g,\bbP) - R(\bbP); \CP \right) \leq R_\ell(g,\bbP) - R_\ell(\bbP). \]
\end{theorem}

{\ctl In general, ensuring uniform calibration of a loss function is much harder than showing Fisher consistency. {\cbl To end this section}, we outline a general strategy to show uniform calibration for generic loss functions, which involves lower-bounding $\delta_\ell$ defined in \eqref{eqn:jpo-uniform-calibration-function}. In Section \ref{sec:U-calib-LS}, we demonstrate this strategy for the squared loss {\cbl $\ell_{\LS}$} for the general class of square-integrable distributions, and then, invoking results from \citet{Steinwart2007}, we show uniform calibration for the class of {\cbl separable} loss functions with respect to the class of {\cbl symmetric} distributions. {\cbl In Section \ref{sec:jpo-uniform-calibration-spo}, we show that for $m=1$, uniform calibration can fail for the $\SPOp$ loss function of \citet{ElmachtoubGrigas2017} even when Fisher consistency is satisfied, and we give a sufficient condition on the class of continuous symmetric distributions that guarantees uniform calibration.}

We first present an alternative form for $\delta_\ell$.
\begin{lemma}\label{lemma:jpo-uniform-calibration-function}
	Consider $\delta_\ell$ defined in \eqref{eqn:jpo-uniform-calibration-function}. We have
	\begin{equation}\label{eqn:jpo-uniform-calibration-alternate}
		\delta_\ell(\epsilon; \CP) = \inf\limits_{x,x' \in X} \inf\limits_{\substack{d : x^*(d) = x\\ \bar{c} : x^*(\bar{c}) = x'}} \inf\limits_{\substack{\bbP \in \CP\\ w \in W\\ \bbE[c \mid w] = \bar{c}}} \left\{ \bbE[\ell(d,c) \mid w] - \min_{d' \in \bbR^m} \bbE[\ell(d',c) \mid w] :~ f(x) - f(x') + \bar{c}^\top (x - x') ~{\cbl \geq}~ \epsilon \right\}.
	\end{equation}
\end{lemma}

We now give a bound on the distance between $d$ and $\bar{c}$ in the second infimum {\cbl in \eqref{eqn:jpo-uniform-calibration-alternate}.}
\begin{lemma}\label{lemma:jpo-uniform-calibration-distance-bound}
	Fix distinct $x,x' \in X$. Let $d$ and $\bar{c}$ be such that $x^*(d) = x$, $x^*(\bar{c}) = x'$. Then
	\[ \|d - \bar{c}\|_2 \geq \frac{\max\{ 0, f(x) - f(x') + \bar{c}^\top (x - x') \}}{\|x - x'\|_2}. \]
\end{lemma}

The strategy to prove $\CP$-calibration of $\ell$ is as follows. First, fixing $x,x'$, notice that if $\bar{c}$ {\cbl and} $d$ are chosen according to the conditions of Lemma \ref{lemma:jpo-uniform-calibration-distance-bound}, together with the condition that $f(x) - f(x') + \bar{c}^\top (x - x') > \epsilon$, then $\|d - \bar{c}\|_2 > \epsilon/\|x - x'\|_2 \geq \epsilon/B_X > 0$ {\cbl holds}, where $B_X$ is the Euclidean diameter of $X$ defined in \eqref{eq:jpo-uniform-diameters}. Then, we want to give a positive lower bound for $\bbE[\ell(d,c) \mid w] - \min_{d' \in \bbR^m} \bbE[\ell(d',c) \mid w]$ over all distributions $\bbP \in \CP$ and $w \in W$ such that $\bbE[c \mid w] = \bar{c}$. {\cbl To this end, we will}  %
exploit the fact that $\argmin_{d' \in \bbR^m} \bbE[\ell(d',c) \mid w]$ is close to $\bbE[c \mid w] = \bar{c}$, and the fact that $\|d - \bar{c}\|_2 > \epsilon/B_X$.}

\subsection{Uniform Calibration of the Squared Loss and Related Loss Functions}\label{sec:U-calib-LS}

We now {\cbl specifically} consider the squared loss function:
\begin{align*}
\ell_{\LS}(d,c) &:= \|d - c\|_2^2\\
\CP &:= \left\{ \bbP :~ \forall w \in W, \ \bbP[\cdot \mid w] \text{ is square integrable, and } \bbE[c \mid w] \in \Conv(C) \right\}.
\end{align*}
Due to the bias-variance decomposition of the squared loss, we can write $\delta_{\ell_{\cbl \LS}}$ entirely as a geometric quantity, without any probabilistic terms.
\begin{lemma}\label{lemma:jpo-delta-simple-squared-loss}
{\cbl Consider the case of the squared loss $\ell_{\LS}$ and $\cP$ as defined above. Then, we}  
have
\[\delta_{\ell_{\LS}}(\epsilon;\CP) = \inf_{x,x' \in X} \inf_{\substack{d \in \bbR^m : x^*(d) = x\\ \bar{c} \in \Conv(C) : x^*(\bar{c}) = x'}} \left\{ \|d-\bar{c}\|_2^2 : f(x) - f(x') + \bar{c}^\top (x - x') ~{\cbl \geq}~ \epsilon \right\}.\]
\end{lemma}

Using Lemmas \ref{lemma:jpo-uniform-calibration-distance-bound} and \ref{lemma:jpo-delta-simple-squared-loss}, we derive $\CP$-uniform calibration of the squared loss {\cbl $\ell_{\LS}$}.
\begin{theorem}\label{thm:jpo-delta-squared-loss}
The squared loss $\ell_{\cbl \LS}$ is $\CP$-uniformly calibrated, with
\[ \delta_{\ell_{\cbl \LS}}(\epsilon; \CP) \geq \frac{\epsilon^2}{B_X^2} > 0 \quad \text{for all $\epsilon > 0$}. \]
\end{theorem}

\begin{corollary}\label{cor:jpo-squared-loss-risk-bound}
{\cbl For} %
 the squared loss {\cbl $\ell_{\LS}$}, we have
\[ \frac{1}{B_X^2}\left( R(g,\bbP) - R(\bbP) \right)^2 \leq R_{\ell_{\LS}}(g,\bbP) - R_{\ell_{\LS}}(\bbP). \]
\end{corollary}

\begin{remark}\label{rem:jpo-squared-loss-implication}
Theorem \ref{thm:jpo-delta-squared-loss} and Corollary \ref{cor:jpo-squared-loss-risk-bound} show that bounding the risk of the squared loss of a predictor $g:W \to \bbR^m$ is enough to bound the true risk. Intriguingly, this holds despite {\cbl the fact that} the squared loss contains %
no information about the optimization problem at hand (i.e., $f$ or $X$). This means that minimization of the true risk can be achieved by training a predictor $g$ without any information on the optimization problem, which is quite counter-intuitive. Furthermore, let $g = (g_1,\ldots,g_m)$ where each $g_j:W \to \bbR$, and observe that
\[ R_{\ell_{\LS}}(g,\bbP) - R_\ell(\bbP) = \sum_{j \in [m]} \left( \bbE[(g_j(w) - c_j)^2] - \inf_{g_j'} \bbE[(g_j'(w) - c_j)^2] \right). \]
Thus, the excess squared loss risk is separable in the coefficients $j \in [m]$, hence we can train individual predictors $g_j:W \to \bbR$ to predict each coefficient $c_j$. Our results state that individual squared error risk bounds are enough to obtain bounds on the true risk $R(g,\bbP)$. {\cbl In particular,} %
invoking Corollary \ref{cor:jpo-squared-loss-risk-bound} gives
\[ R(g,\bbP) - R(\bbP) \leq B_X \sqrt{\sum_{j \in [m]} \left( \bbE[(g_j(w) - c_j)^2] - \inf_{g_j'} \bbE[(g_j'(w) - c_j)^2] \right)}. \]
Again, this is quite counter-intuitive, since we know that {\cbl a small change in only} %
one coefficient of $d$ can change the optimal solution $x^*(d)$.
\epr
\end{remark}

Remark \ref{rem:jpo-squared-loss-implication} states that squared error risk bounds on individual coefficients $j \in [m]$ are enough to bound the true optimality gap risk, which essentially states that one-dimensional least squares regression on each coefficient $j \in [m]$ is sufficient for end-to-end prediction and optimization. Several other loss functions have been utilized in regression, due to their superior finite-sample performance. For example, the absolute deviation loss from Example \ref{ex:jpo-abs-dev-loss} or the Huber loss have been used for heavy-tailed data due to their reduced sensitivity to outliers. \citet[Section 4.3]{Steinwart2007} studies the use of alternate loss functions in regression, and their risk relationships to the squared loss risk. By invoking these results, we can correspondingly obtain bounds on the true risk. More precisely, we have the following result.
\begin{lemma}\label{lemma:jpo-uniform-calibration-decompose}
For each $j \in [m]$, let $\ell_j : \bbR \times \bbR \to \bbR$ be a loss function such that there exists a non-decreasing function $\eta_j:(0,\infty) \to (0,\infty)$ that satisfies
\begin{equation}\label{eqn:jpo-squared-loss-one-dim}
\bbE[(g_j(w) - c_j)^2] - \inf_{g_j'} \bbE[(g_j'(w) - c_j)^2] \leq \eta_j\left( \bbE[\ell_j(g_j(w),c_j)] - \inf_{g_j'} \bbE[\ell_j(g_j(w),c_j)] \right)
\end{equation}
for any $g_j:W \to \bbR$ and $\bbP \in \CP$. Then, denoting $g = (g_1,\ldots,g_m)$,
\[ R(g,\bbP) - R(\bbP) \leq B_X \sqrt{\sum_{j \in [m]} \eta_j\left( \bbE[\ell_j(g_j(w),c_j)] - \inf_{g_j'} \bbE[\ell_j(g_j(w),c_j)] \right)}. \]
\end{lemma}

Thus, when we use a separable loss function $\ell(d,c) = \sum_{j \in [m]} \ell_j(d_j,c_j)$ to train a predictor $g$, we can obtain true risk bounds by deriving bounds on each $\bbE[\ell_j(g_j(w),c_j)] - \inf_{g_j'} \bbE[\ell_j(g_j(w),c_j)]$. Conditions for the existence of the functions $\eta_j$ are given by results from \citet{Steinwart2007}. To obtain them, we need to restrict the class of distributions. Precisely, we define $\CP_{\sym}$ to be the class of square integrable distributions such that for all $w \in W$ and $j \in [m]$, $\bbP[c_j \mid w]$ is a symmetric distribution, i.e., $c_j - \bbE[c_j \mid w]$ has the same conditional distribution as $\bbE[c_j \mid w] - c_j$. 
\begin{theorem}[{\citet[Theorems 4.19, 4.20(ii)]{Steinwart2007}}]\label{thm:jpo-squared-loss-calibrate}
Fix any $j \in [m]$. Let $\ell_j(d_j,c_j) = \psi_j(d_j - c_j)$ where $\psi_j: \bbR \to [0,\infty)$ is symmetric, i.e., $\psi_j(r) = \psi_j(-r)$, and uniformly convex, i.e., there exists some non-decreasing $\delta_j:[0,\infty) \to [0,\infty)$ with $\eta_j(0) = 0$ such that for all $\alpha \in [0,1]$ and $r,r' \in \bbR$,
\[ \alpha \psi_j(r) + (1-\alpha) \psi_j(r') - \psi_j(\alpha r + (1-\alpha)r') \geq \alpha(1-\alpha) \delta_j(|r - r'|^2). \]
Then, for any $g_j:W \to \bbR$ and $\bbP \in \CP_{\sym}$, we have
\[ \frac{1}{4} \delta_j^{**}\left( \bbE[(g_j(w) - c_j)^2] - \inf_{g_j'} \bbE[(g_j'(w) - c_j)^2] \right) \leq \bbE[\psi_j(g_j(w) - c_j)] - \inf_{g_j'} \bbE[\psi_j(g_j(w) - c_j)]. \]
\end{theorem}
While the proof {\cbl of Theorem~\ref{thm:jpo-squared-loss-calibrate}} can be found in the relevant sections of \citet{Steinwart2007}, we give a more concise version in \ref{sec:jpo-proofs-non-asymptotic}.

\subsection{Uniform Calibration of the  {\cbl $\SPOp$} Loss in Example \ref{ex:jpo-spo-one-dim}}\label{sec:jpo-uniform-calibration-spo}

Recall {\cbl our Example \ref{ex:jpo-spo-loss} that studied} the {\cbl $\SPOp$} loss \eqref{eqn:jpo-spo-loss} introduced in \citet[Definition 3]{ElmachtoubGrigas2017}. {\cbl It was shown in} \citet[Theorem 1]{ElmachtoubGrigas2017} that this loss is Fisher consistent, hence by Theorem \ref{thm:jpo-FisherConsistent-calibration} it is $\bbP$-calibrated {\cbl whenever} $\bbP[c \mid w]$ is centrally symmetric and continuous for all $w \in W$. {\cbl On the other hand, the uniform calibration of the $\SPOp$ loss~\eqref{eqn:jpo-spo-loss} has not yet been studied. }%
In this section, we {\cbl examine its} uniform calibration for the special one-dimensional case $m=1$, i.e., Example \ref{ex:jpo-spo-one-dim}; to our knowledge, the general $m$ case remains open.

Recall Example \ref{ex:jpo-spo-one-dim} has $m = 1$, $f = 0$, $X = [-1/2,1/2]$, and we will take $C = \bbR$. In this case, {\cbl recall that} the loss function \eqref{eqn:jpo-spo-loss} becomes
\[ \ell_{\SPOp}(d,c) := \frac{1}{2} \left( |2d-c| -2d \sign(c) + |c| \right). \]
{\cbl For this loss function, \citet{ElmachtoubGrigas2017} studied a particular class of probability distributions that are} symmetric and continuous over $\bbR$. 
{\cbl Recall that a} continuous distribution is one such that the probability density function (w.r.t.\@ Lebesgue measure) is positive over all of $\bbR$. For simplicity, we consider {\cbl the same class of} symmetric, continuous distributions over $\bbR$, {\cbl i.e.,} %
\[ \CP_{\cont,\sym} := \left\{ \bbP : \ \forall w \in W,\ \bbP[c \mid w] \text{ is continuous and symmetric} \right\}. \]
{\cbl For this class of distributions, in} \citet{ElmachtoubGrigas2017} {\cbl it was shown} that the (conditional) mean is the unique minimizer of $\min_{d' \in \bbR} \bbE[\ell(d',c) \mid w]$.
\begin{lemma}[{\citet[Theorem 1]{ElmachtoubGrigas2017}}]\label{lemma:jpo-spo-loss-Fisher}
Let $\bbP \in \CP_{\cont,\sym}$. Then for any $w \in W$, the unique minimizer of $\min_{d' \in \bbR} \bbE[\ell_{\SPOp}(d',c) \mid w]$ is $d^* = \bbE[c \mid w]$.
\end{lemma}

Using Lemmas \ref{lemma:jpo-uniform-calibration-function} and \ref{lemma:jpo-spo-loss-Fisher}, and {\cbl noting} %
\[ x^*(d) = \begin{cases}
-1/2, & d > 0\\
0, & d = 0\\
1/2, & d < 0,
\end{cases} \]
we have
\[ \delta_{\ell_{\SPOp}}(\epsilon;\CP_{\cont,\sym}) = \inf_{w \in W} \inf_{\substack{d,\bar{c} \in \bbR : d\bar{c} < 0}} \inf_{\substack{\bbP \in \CP_{\cont,\sym}\\ \bbE[c \mid w] = \bar{c}}} \left\{ \bbE[\ell_{\SPOp}(d,c) \mid w] - \bbE[\ell_{\SPOp}(\bar{c},c) \mid w] :~ |\bar{c}| > \epsilon \right\}. \]
Fixing $w \in W$, assume that $\bbE[c \mid w] = \bar{c} > 0$, hence $d < 0 < \epsilon < \bar{c}$. Since {\cbl the function} $\bbE[\ell_{\SPOp}(d,c)] = \frac{1}{2} \left( \bbE[|2d-c| \mid w] - 2d (\bbP[c > 0 \mid w] - \bbP[c < 0 \mid w]) + \bbE[|c| \mid w] \right)$ is convex in $d$ and hence continuous, when restricting $d < 0$, the closest $\bbE[\ell_{\SPOp}(d,c) \mid w]$ can get to the minimizer $\bbE[\ell_{\SPOp}(\bar{c},c) \mid w]$ is at $d = 0$, i.e., $\bbE[\ell_{\SPOp}(0,c) \mid w] = \bbE[|c| \mid w]$. A similar argument holds for $\bar{c} < 0$. Therefore, 
\begin{align*}
\delta_{\ell_{\SPOp}}(\epsilon;\CP_{\cont,\sym}) &= \inf_{w \in W} \inf_{\substack{|\bar{c}| > \epsilon\\ \bbP \in \CP_{\cont,\sym}\\ \bbE[c \mid w] = \bar{c}}} \left\{ \bbE[\ell_{\SPOp}(0,c) \mid w] - \bbE[\ell_{\SPOp}(\bar{c},c) \mid w] \right\}\\
&= \inf_{w \in W} \inf_{\substack{|\bar{c}| > \epsilon\\ \bbP \in \CP_{\cont,\sym}\\ \bbE[c \mid w] = \bar{c}}} \left\{ \bbE[|c| \mid w] - \frac{1}{2} \bigg( \bbE[|2\bar{c}-c| \mid w] + 2 \bar{c} \left( \bbP[c > 0 \mid w] - \bbP[c < 0 \mid w] \right) + \bbE[|c| \mid w] \bigg) \right\}\\
&= \inf_{w \in W} \inf_{\substack{\bbP \in \CP_{\cont,\sym}\\ |\bbE[c \mid w]| > \epsilon}} \left\{ \bbE[c \mid w] \left( \bbP[c > 0 \mid w] - \bbP[c < 0 \mid w] \right) \right\},
\end{align*}
where the third equality follows because $\bbP[c \mid w]$ is symmetric, so $2 \bar{c} - c$ has the same conditional distribution as $c$, thus $\bbE[|2\bar{c}-c| \mid w] = \bbE[|c| \mid w]$. Unfortunately, we can show that $\delta_\ell(\epsilon;\CP_{\cont,\sym}) = 0$ for all $\epsilon > 0$. {\cbl This is due to the following result.}
\begin{proposition}\label{prop:jpo-spo-one-dim-non-calibrated}
For any $\epsilon > 0$, {we can construct} a sequence of symmetric, continuous distributions $\left\{ \bbP^{(k)} \right\}_{k \in \bbN}$ on $\bbR$ with $|\bbE^{(k)}[c]| \geq \epsilon$ such that $\bbE^{(k)}[c] \left( \bbP^{(k)}[c > 0] - \bbP^{(k)}[c < 0] \right) \to 0$. {\cbl Therefore, by Lemma \ref{lemma:jpo-uniform-calibration}, $\ell_{\SPOp}$ is not $\CP_{\cont,\sym}$-calibrated even in the restricted $m=1$ setting.}
\end{proposition}

{\cbl In contrast to this, we close this section by establishing a uniform calibration result for $\ell_{\SPOp}$ for the case of the more restrictive class of continuous and symmetric distributions with uniformly bounded margin $\left| \bbP[c > 0 \mid w] - \bbP[c < 0 \mid w] \right|$.} %
\begin{proposition}\label{prop:jpo-spo-one-dim-calibrated}
For $\alpha > 0$, let
\[\CP_{\cont,\sym,\alpha} := \left\{ \bbP : \ \forall w \in W,~ \begin{aligned}
&\bbP[c \mid w] \text{ is continuous and symmetric}\\
&\left| \bbP[c > 0 \mid w] - \bbP[c < 0 \mid w] \right| \geq \alpha
\end{aligned} \right\}.\]
Then, $\ell_{\SPOp}$ is $\CP_{\cont,\sym,\alpha}$-calibrated, and we have
\[ R(g,\bbP) - R(\bbP) \leq \frac{1}{\alpha} \left( R_{\ell_{\SPOp}}(g,\bbP) - R_{\ell_{\SPOp}}(\bbP)\right). \]
\end{proposition}

\section{Computational Study}\label{sec:jpo-computational}

{\cbl In this section, we conduct a computational study in order to investigate the effect of consistency in end-to-end prediction and optimization frameworks, and the effect of using a loss function that takes into account the optimization problem information.} 
{\cbl For this purpose, we examine the squared loss $\ell_{\LS}$ and the $\SPOp$ loss $\ell_{\SPOp}$ in our experiments. Recall that}
the squared loss {\cbl $\ell_{\LS}$} does not take into account any information about the optimization problem, e.g., $f$ or $X$, {\cbl yet in} Theorem \ref{thm:jpo-delta-squared-loss} and Corollary \ref{cor:jpo-squared-loss-risk-bound} {\cbl we} provided true risk bounds in terms of the surrogate squared loss risk bounds. {\cbl In contrast,} \citet{ElmachtoubGrigas2017} {\cbl proposed} the {\cbl $\SPOp$ loss} function {\cbl $\ell_{\SPOp}$} \eqref{eqn:jpo-spo-loss}, which incorporates information about the optimization problem, and is known to be Fisher consistent with respect to certain distributions (see Example \ref{ex:jpo-spo-loss}) {\cbl but has weaker calibration properties than $\ell_{\LS}$ (see Section \ref{sec:jpo-uniform-calibration-spo})}.

{\cbl
Recall also that the squared loss $\ell_{\LS}$ and the $\SPOp$ loss $\ell_{\SPOp}$ are defined as follows:
\begin{align*}
\ell_{\LS}(d,c) &:= \|d-c\|_2^2\\
\ell_{\SPOp}(d,c) &:= f(x^*(c)) + (2d-c)^\top x^*(c) - \min_{x \in X} \left\{ f(x) + (2d-c)^\top x \right\} = L(c,2d-c).
\end{align*}
Note that \citet[Section 3.2]{ElmachtoubGrigas2017} originally defined the $\SPOp$ loss for linear objectives $c^\top x$ only, with $f=0$. However, the above definition is a straightforward extension of their derivation for the objective $f(x) + c^\top x$.
}

{\cbl We investigate three problem classes. First, we examine portfolio optimization using real-world data, where consistency is not known a priori. Second, we examine the fractional knapsack problem on simulated data, where some chosen parameters control the degree of non-linearity of the underlying data model, and thereby the consistency of certain loss functions. Third, we examine multiclass classification on simulated data, where the $\SPOp$ loss is provably inconsistent (see Example \ref{ex:jpo-multiclass}), but the squared loss is consistent.
	
In all problem classes, we compare linear predictors $w \mapsto g(w) := Vw$ where $V$ is obtained by solving the empirical risk minimization problem
\begin{equation}\label{eq:jpo-ERM-linear}
\min_{V \in \bbR^{m \times k}} \frac{1}{n} \sum_{i \in [n]} \ell(Vw_i,c_i)
\end{equation}
for different loss functions $\ell$ on the same historical data $\left\{ (w_i,c_i) \right\}_{i \in [n]}$. 
}

{\cbl
Our results suggest the following key managerial insights:
\begin{itemize}
	\item On portfolio instances constructed from real data, there is no significant difference between using $\ell_{\LS}$ loss and the more high-powered $\SPOp$ loss which takes into account optimization problem information.
	
	\item Recall that we have shown theoretically that on the multiclass classification problem, the $\SPOp$ loss is provably inconsistent (Example \ref{ex:jpo-multiclass}). Moreover, our numerical results on these problem instances show that %
	the performance of $\SPOp$ loss (expectedly) deteriorates. This highlights the importance of ensuring a consistent loss function is used in practice whenever possible.
	
	\item On the fractional knapsack instances, we observe that constructing a close approximation of the true loss by regularization, while conceptually reasonable, does not provide good empirical results due to the considerable increase in computational effort required to find the corresponding estimator. Therefore, computational efficiency plays an important role in the prediction and optimization context.
	
	\item Our experiments on the fractional knapsack instances also highlight that there are further properties besides consistency and calibration that can be investigated, such as robustness to model misspecification, where $\SPOp$ has an advantage.
\end{itemize}
}

\begin{figure}[tb]
\centering
\includegraphics[page=2,scale=1]{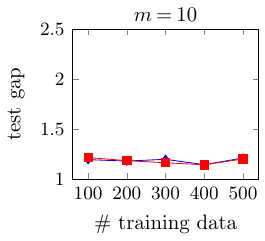}
\includegraphics[page=3,scale=1]{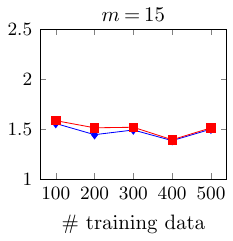}
\includegraphics[page=4,scale=1]{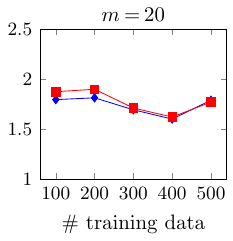}
\includegraphics[page=5,scale=1]{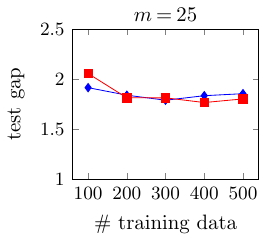}
\includegraphics[page=6,scale=1]{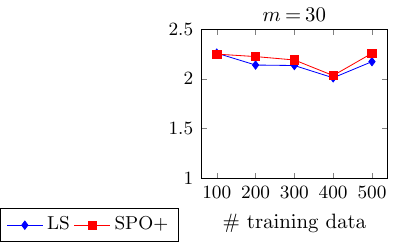}
\caption{Median test optimality gap for different $m$ for portfolio optimization.}
\label{fig:jpo-portfolio-test-gap}
\end{figure}

{\cbl 
\subsection{Mean-variance portfolio optimization}\label{sec:portfolio}

The mean-variance portfolio optimization problem can be expressed as the following constrained quadratic optimization problem
\[ \min_{x \in X} \left\{ f(x) - c^\top x \right\}, \quad \text{where}\quad  f(x) = \frac{1}{2} x^\top Q x, \quad X := \left\{ x \in \bbR^m :~ p^\top x = b,~ x \geq 0 \right\}, \] 
and $Q \succ 0$ is positive definite matrix. This problem arises from portfolio optimization: $x$ denotes a vector of weights for each asset which specifies what proportion of our wealth to investigate in each one, the random vector $c$ represents returns of each stock, the quadratic term $f(x) = \frac{1}{2} x^\top Q x$ represents the risk of the portfolio (usually its variance), and a wealth constraint is imposed with $p = \bm{1}$ and $b = 1$.

In our study, we assume that $c$ is uncertain but $Q$ is fixed and known. 
We follow the common hypothesis in portfolio optimization that the expected cost vector can be described via a linear model $\bbE[c \mid \tilde{w}] = \tilde{b} + \tilde{V}\tilde{w}$, where $\tilde{w}$ are market factors (see \citep{FamaFrench1992}). In this setting, $\tilde{b}$ is the mean vector and $\tilde{V}$ is called the `factor loading matrix.' The goal in this problem is to estimate both $\tilde{b}$ and $\tilde{V}$. To simplify notation, we append a $1$ to each feature vector and denote $w = (\tilde{w},1)$. Similarly, we add $\tilde{b}$ as a column to $\tilde{V}$, and denote $V = (\tilde{V},\tilde{b})$. Thus, our model is $\bbE[c \mid w] = Vw$, and we aim to estimate $V$. We do this by again minimizing \eqref{eq:jpo-ERM-linear} where we take $\ell$ to be  $\ell_{\LS}$ or $\ell_{\SPOp}$. Note that for objectives of type $f(x) - c^\top x$, the $\SPOp$ loss is
\begin{align*}
	\ell_{\SPOp}(d,c) &= L(c,2d-c) = f(x^*(c)) - (2d-c)^\top x^*(c) - \min_{x \in X} \left\{ f(x) - (2d-c)^\top x \right\}.
\end{align*}

Usually, in portfolio optimization, we are permitted to have entries of $x$ negative, which means we short-sell some assets. We show that if we redefine the domain to be $X := \left\{ x \in \bbR^m :~ p^\top x = b \right\}$ without the non-negativity constraints, the $\SPOp$ loss and the true loss are equivalent.
\begin{proposition}\label{prop:jpo-portfolio-SPOp-longshort}
	Let $X := \left\{ x \in \bbR^m :~ p^\top x = b \right\}$ and $A := Q^{-1} - \frac{1}{p^\top Q^{-1} p} Q^{-1} p (Q^{-1} p)^\top$. Then, for any $d$, the optimal solution to $\min_{x \in X} \left\{ \frac{1}{2} x^\top Q x - d^\top x \right\}$ is
	\[ x^*(d) = A d + \frac{b}{p^\top Q^{-1} p} Q^{-1} p. \]
	Furthermore,
	\[ L(d,c) = \frac{1}{2} x^*(d)^\top Q x^*(d) - c^\top x^*(d) - \min_{x \in X} \left\{ \frac{1}{2} x^\top Q x - c^\top x \right\} = \frac{1}{2} (d-c)^\top A (d-c). \]
	Consequently,
	\begin{align*}
		\ell_{\SPOp}(d,c) &= L(c,2d-c) = 2 (c-d)^\top A (c-d) = 4 L(d,c).
	\end{align*}	
\end{proposition}

When we consider linear predictors $w \mapsto Vw$, we can show that the least squares loss $\ell_{\LS}(d,c) = \frac{1}{2} \|d-c\|_2^2$ also optimizes the true loss. More precisely, given data $\{(w_i,c_i) : i \in [n]\}$, a solution to $\frac{1}{n} \sum_{i \in [n]} L(V w_i, c_i)$ can be obtained by minimizing $\frac{1}{n} \sum_{i \in [n]} \ell_{\LS}(V w_i, c_i)$.
\begin{proposition}\label{prop:jpo-portfolio-least-squares}
	Given a matrix $A \succeq 0$ and random variables $(w,c) \sim \bbP$ such that $\bbE[ww^\top]$ is invertible, we have
	\[ \argmin_{V} \bbE\left[ \frac{1}{2} (Vw - c)^\top A (Vw - c) \right] = \argmin_{V} \bbE\left[ \frac{1}{2} \|Vw - c\|_2^2 \right] + \left\{ \tilde{V} :~ A \tilde{V} = \bm{0} \right\}. \]
	Consequently, when  $X = \left\{ x \in \bbR^m : p^\top x = b \right\}$ and $A = Q^{-1} - \frac{1}{p^\top Q^{-1} p} Q^{-1} p (Q^{-1} p)^\top$, the minimizers of $\bbE[\ell_{\LS}(Vw,c)]$ are also minimizers of $\bbE[L(Vw,c)]$.
\end{proposition}
Proofs of Propositions \ref{prop:jpo-portfolio-SPOp-longshort} and  \ref{prop:jpo-portfolio-least-squares} are in Section \ref{sec:jpo-proofs-computational}.

For this reason, in our numerical study we henceforth impose non-negativity constraints on our decision variables $X := \left\{ x \in \bbR^m : x \geq 0, \ p^\top x = b \right\}$. We generate instances from data on stocks that remained in the S\&P 500 index for all 1258 trading days between January 1, 2003 and December 31, 2007. We also collected data on the three Fama-French factors for these trading days, these are our feature vectors, with a $1$ appended, so $k=4$.

We consider $m \in \{10,15,\ldots,30\}$, and for each $m$, we generate 100 random instances by choosing $m$ random stocks. For each instance, we collect $n \in \{100,\ldots,500\}$ consecutive days of stock returns for the set of chosen stocks; stock returns for a particular day are recorded as the percentage increase/decrease of that day's price from the previous day's price. The matrix $Q$ is the $m \times m$ sample covariance matrix of the stock returns computed from the $n$ training days. We then estimate $V$ from the $n$ days of stock returns data via optimizing the least squares loss and the $\SPOp$ loss. We evaluate the performance of our estimated $V$ on the next $N=10$ days after the $n$-day window in the training data, by first taking the factor data $w$ for each test day, computing $Vw$, using that to compute a portfolio $x^*(Vw)$, then computing the objective of that portfolio on the actual $f(x^*(Vw)) - c^\top x^*(Vw)$ for that day.
We report the median optimality gap $L(Vw,c) = f(x^*(Vw)) - c^\top x^*(Vw) - \left( f(x^*(c)) - c^\top x^*(c) \right)$ (so lower is better) in Figure \ref{fig:jpo-portfolio-test-gap}, which shows little difference between using the $\SPOp$ loss and least squares on this class of problems with real data.
}

\subsection{Fractional knapsack problem}\label{sec:knapsack}

{\cbl In the case of} %
fractional knapsack linear programs, %
{\cbl we have}
\begin{equation}\label{eqn:jpo-fractional-knapsack-def}
\max_{x \in X} d^\top x, \quad \text{where} \quad X := \left\{ x \in [0,1]^m : p^\top x \leq B \right\}, \quad \text{and} \quad f(x) = 0.
\end{equation}
{\cbl Here,} $p \in \bbR^m$ is some fixed positive vector, and $B > 0$ is the capacity of the knapsack. {\cbl As before, we test $\ell_{\LS}$ and $\ell_{\SPOp}$. Note that due to the max-type optimization problem, the $\SPOp$ loss becomes
\[ \ell_{\SPOp}(d,c) = \max_{x \in X} (2d-c)^\top x - (2d-c)^\top x^*(c) = L(c,2d-c). \]
For this problem class, we also test an additional loss function
\begin{align*}
\ell_{\reg,\lambda}(d,c) &:= c^\top x^*(c) - c^\top x_{\lambda}^*(d), \quad x_{\lambda}^*(d) := \argmax_{x \in X} \left\{ d^\top x - \frac{\lambda}{2} \|x\|_2^2 \right\}.
\end{align*}
Note that the loss function $\ell_{\reg,\lambda}$ is  nothing but the \emph{exact} optimality gap evaluated at the \emph{unique} solution to the regularized knapsack problem with  the objective function $d^\top x - \frac{\lambda}{2} \|x\|_2^2$  that includes a regularization term}. We consider the regularized objective due to the fact that %
the set of optimal solutions $X^*(d)$ for the unregularized problem {\cbl does not admit a simple model}. By adding a regularizer, however, we can show that $\ell_{\reg,\lambda}(d,c)$ is mixed-integer linear representable.
{\cbl 
\begin{proposition}\label{prop:jpo-knapsack-regularized-mip}
For fixed $c$ and $\lambda$, The set $\{ (d,t) : \ell_{\reg,\lambda}(d,c) \leq t \}$ admits a mixed-integer linear representation. Consequently, the empirical risk minimization problem \eqref{eq:jpo-ERM-linear} with $\ell = \ell_{\reg,\lambda}$ can be formulated as a mixed-integer linear program.
\end{proposition}
The proof of Proposition~\ref{prop:jpo-knapsack-regularized-mip} is rather standard;} thus we give the details in Section \ref{sec:jpo-proofs-computational}.

We {\cbl generate and test} knapsack instances with $m=10$ {\cbl with synthetic data} as follows. Each item weight $p_j$ is a random integer between $1$ and $1000$. Then, $B$ is a random integer between $l$ and $u$, where $l = \max_{j \in [m]} p_j$, $u = (r l/\bm{1}^\top p + 1 - l/\bm{1}^\top p) \bm{1}^\top p$, where $r$ is uniformly distributed on $[0,1]$. For $m=10$, $k=5$, we generate 30 knapsack instances in this way, each paired with a randomly chosen {\cbl ground truth coefficient matrix $V_0 \in \bbR^{m \times k}$}. To generate data from $V_0$, we use a similar scheme to that of \citet{ElmachtoubGrigas2017}. The feature support set is $W := [-1,1]^k$, and each $w_i$ is drawn uniformly at random from $W$, except that the last entry $w_{ik}=1$ always (in this way we can model a constant term in our predictor). Then, given hyper-parameters $\delta \geq 1,\epsilon \in (0,1)$, each $c_i$ is generated as
\[ c_{ij} := \tilde{\epsilon}_{ij} (v_{0,j}^\top w)^\delta + \eta_{ij}, \quad j \in [m] \]
where $\tilde{\epsilon}_{ij}$ is uniformly distributed on $[1-\epsilon,1+\epsilon]$ and $2\eta_{ij}+1$ is an exponential random variable with scale parameter $\lambda = 1$ (thus $\eta_{ij}$ has zero mean). Note that the exponentiation by $\delta$ is entry-wise, and that when $\delta = 1$ this means we have a linear model with random noise. We test $\delta = 1,3,5$ and $\epsilon = 0.1$ for each instance. We consider datasets of size $n=100,200,300,400,500$ generated in this way. We trained $\ell_{\reg,\lambda}$ with $\lambda = 0.01$.

\begin{figure}[tb]
	\centering
	\includegraphics[page=7,scale=1]{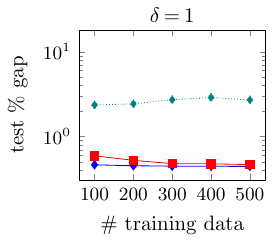}
	\includegraphics[page=8,scale=1]{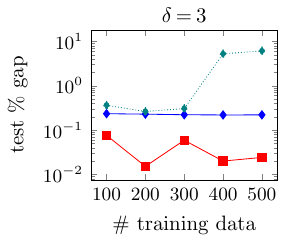}
	\includegraphics[page=9,scale=1]{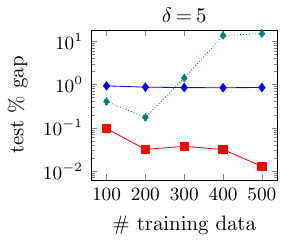}
	\includegraphics[page=10,scale=1]{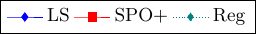}
	\caption{Average test relative optimality gap for different $\delta$ and $\epsilon = 0.2$ for the  continuous knapsack problem.}
	\label{fig:jpo-average-test-gap}
\end{figure}

To test our predictors, we generate $10,000$ points from the same distribution for each hyper-parameter setting and $V_0$, and evaluate the average optimality gap using $L$ on the test set for our predictors. Our results are shown in Figure \ref{fig:jpo-average-test-gap} where we measure the average percentage optimality gap $\tilde{L}(d,c) = c^\top (x^*(c) - x^*(d))/(c^\top x^*(c))$ (so lower is better).

First, it is clear that $\ell_{\reg,\lambda}$ has poorer performance than $\ell_{\LS}$ and $\ell_{\SPOp}$. We attribute this to the fact that very few problems were {\cbl solved} to optimality within the five minute time limit. {\cbl Therefore, this brings up the insight that  despite $\ell_{\reg,\lambda}$ being a close approximation to the true loss $L$ on paper, computational considerations must be taken into account during training. Second, notice that for higher values of $\delta$ (i.e., as the true model becomes more non-linear), $\ell_{\SPOp}$ outperforms $\ell_{\LS}$, which points to a `robustness to prediction model misspecification' property that $\ell_{\SPOp}$ might satisfy, and suggests that taking into account optimization information may increase performance under model misspecification. This phenomenon of $\ell_{\SPOp}$ is currently unexplained by the theoretical results, and is an interesting direction for future research.
}

{\cbl
\subsection{Multiclass classification}\label{sec:num-multiclass}

In our last class of examples, we consider the setting of multiclass classification from Example~\ref{ex:jpo-multiclass} with $C = \{c_j := \bm{1}_m - e_j : j \in [m]\} \subset \bbR^m$, $X = \Conv\left\{ e_j : j \in [m] \right\} \subset \bbR^m$ and $f(x) = 0$ for all $x \in X$. Recall that $e_j \in \bbR^m$ denotes the $j$th standard basis vector for $j \in [m]$.  The $\SPOp$ loss for this problem class is given by 
\begin{align*}
	\ell_{\SPOp}(d,c_j) &= (2d-c_j)^\top e_j - \min_{x \in X} (2d - c_j)^\top x = 2d_j - \min_{x \in X} \left( 2d_j x_j + \sum_{j' \in [m], j' \neq j} (2d_{j'}-1) x_j \right).
\end{align*}
A lifted representation of the $\SPOp$ loss is given by the following proposition.
\begin{proposition}\label{prop:multiclass-SPOp-loss}
	For fixed $c=c_j$, the set $\left\{ (d,t) :~ \ell_{\SPOp}(d,c_j) \leq t \right\}$
	has a lifted representation
	\[ \left\{ (d,t,\gamma) :~ \begin{aligned}
		&2d_j - \gamma \leq t\\
		&\gamma \leq 2d_j\\
		&\gamma \leq 2d_{j'} - 1, \ j' \in [m] \setminus \{j\}
	\end{aligned} \right\}. \]
\end{proposition}

Recall that in Example~\ref{ex:jpo-multiclass} we have shown $\SPOp$ to be inconsistent for this problem theoretically. We next numerically compare the performance of $\ell_{\LS}$ with $\ell_{\SPOp}$ to investigate the effects of using an inconsistent loss function versus a consistent one. We use simulated data generated under the following multinomial logit model with parameters $v_1,\ldots,v_m \in \bbR^k$:
\[ \bbP[c=c_j \mid w] = \frac{\exp(-v_j^\top w)}{\sum_{j' \in [m]} \exp(-v_{j'}^\top w)} \quad j \in [m]. \]
Under this model, given $w$, choosing the most likely class is equivalent to choosing the index $j$ which gives the smallest $v_j^\top w$. We generate $w$ uniformly at random from the unit cube $[0,1]^k$. We fix $k=4$ and $m=4$. We do 100 repetitions of the following:
\begin{itemize}
	\item Generate a true coefficient matrix $V^{\text{true}} \in \bbR^{m \times k}$ where each entry is distributed as a standard normal random variable.
	\item Generate test features $\{w_i^{\text{test}}\}_{i \in [N]}$ where $N=100,000$. and compute the true probabilities $\{ p_j(w_i^{\text{test}}):=\bbP[c=c_j \mid w_i^{\text{test}}] \}_{j \in [m], i \in [N]}$ using the true parameters $V^{\text{true}}$.
	\item For each $n \in \{100,200,\ldots,1000\}$:
	\begin{itemize}
		\item Generate training data $\{w_i^{\text{train}},c_i^{\text{train}}\}_{i \in [n]}$ according to the true model.
		\item Estimate the parameters using the two proposed methods to obtain $V_{\LS}, V_{\SPOp}$.
		\item Use the test data to evaluate estimated parameters $\hat{V}$ by computing
		\[ \frac{1}{N} \sum_{i \in [N]} \left( 1 - \frac{1}{\left| \argmin_{j' \in [m]} \hat{v}_j^\top w_i^{\text{test}} \right|} \sum_{j \in \argmin_{j' \in [m]} \hat{v}_j^\top w_i^{\text{test}}} p_j\left( w_i^{\text{test}} \right) \right). \]
	\end{itemize}
\end{itemize}

Note that the term in the outer summand is simply $\bbE\left[ L(\hat{V} w_i^{\text{test}}, c) \mid w_i^{\text{test}} \right]$, the expected true loss of plugging in the vector $\hat{V} w_i^{\text{test}}$ into the optimization problem, and if it has a non-unique minimizer then one is chosen at random from the set of minimizers. We can estimate the best possible loss if we had true knowledge of the distribution, i.e., the Bayes loss, as
\[ L_{\Bayes} = \frac{1}{N} \sum_{i \in [N]} \left( 1 - \max_{j \in [m]} p_j(w_i^{\text{test}}) \right). \]

\begin{figure}[h!t]
	\centering
	\includegraphics[page=11,scale=1]{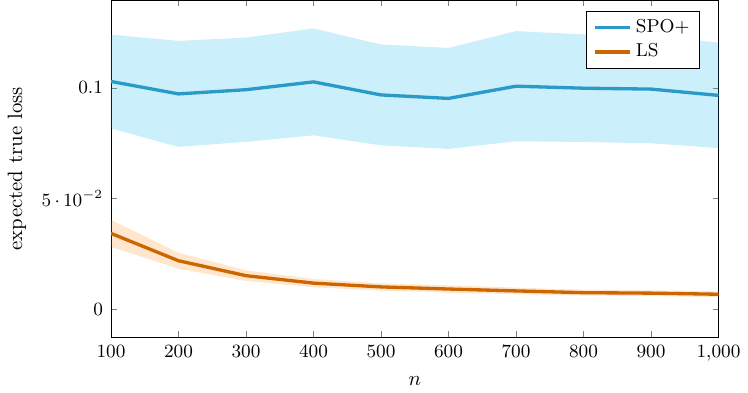}
	\caption{Mean (line) and a two standard deviation range (shaded region) of expected true loss across 100 runs for multiclass classification with $m=k=4$.}
	\label{fig:jpo-multiclass-test}
\end{figure}

In Figure \ref{fig:jpo-multiclass-test} we plot the mean and a two standard deviation band for the gap of the true losses for each predictor relative to the Bayes loss across 100 runs, that is we  plot statistics for the following quantity:
\[\frac{\bbE[L(\hat{V} w, c)] - L_{\Bayes}}{L_{\Bayes}}.\]
It is clear from Figure \ref{fig:jpo-multiclass-test} that the $\SPOp$ loss performs noticeably worse than the least squares loss. This observation is perhaps expected from our theoretical findings since we established that the $\SPOp$ loss is inconsistent for this problem class. However, note that this performance difference between $\LS$ and  $\SPOp$ losses is still interesting because the true (conditional) expected cost vector
\[ \bbE[c \mid w] = \left\{ 1 - \frac{\exp(-v_j^\top w)}{\sum_{j' \in [m]} \exp(-v_{j'}^\top w)} \right\}_{j \in [m]} \]
is a highly non-linear function of $w$ and restricting it to a linear model, such as the case of $V_{\LS} w$, may prevent us from learning the true functional form of $\bbE[c \mid w]$. A potential reason for the superior performance of the least squares loss is that it is consistent. In particular, for a given $w$, even though $V_{\LS} w$ may not exactly be $\bbE[c \mid w]$, the minimal entry may still coincide. On the other hand, we showed in Example \ref{ex:jpo-multiclass} in Section \ref{sec:jpo-proofs-FisherConsistent} that the true minimizer of $\bbE[\ell_{\SPOp}(d,c) \mid w]$ is a constant vector, which we know will not give us the correct minimal entry of $\bbE[c \mid w]$. Our experiments thus highlight an important insight: consistency of a loss function matters more than whether the loss function takes into account optimization problem information. In particular, despite the fact that the $\SPOp$ loss takes into account information from the optimization problem, its inconsistency for this problem class resulted in poor performance.

Notice also that there is no downward trend in the expected true loss of $\SPOp$ as $n$ increases. This is because of its inconsistency. In fact, a closer look at the estimated $V_{\SPOp}$ reveals that it often estimates a zero matrix, which predicts the zero vector $V_{\SPOp} w = \bm{0}$. This is consistent with the theoretical analysis of Example \ref{ex:jpo-multiclass} in Section \ref{sec:jpo-proofs-FisherConsistent}, where it is shows that constant vectors $d$ minimize $\bbE[\ell_{\SPOp}(d,c)]$ when $\max_{j \in [m]} p_j < 1/2$.
}

\section{Conclusion}\label{sec:conclusion}

In this paper, we explored risk guarantees for end-to-end prediction and optimization processes, which are prevalent in practice.
We showed that the true non-convex optimality gap risk can be minimized via minimizing the surrogate risk as long as the surrogate loss function is appropriately calibrated, and provided precise relationships between the two risks under these assumptions. {\cbl We provided an equivalence result (Theorem \ref{thm:jpo-FisherConsistent-calibration}) that allows us to easily check the weaker $\bbP$-calibration condition via Fisher consistency, and used it to explore calibration conditions for certain loss functions in Section \ref{sec:jpo-FisherConsistent-surrogate-risk}. We also examined a stronger notion of uniform calibration for the least squares $\ell_{\LS}$ and $\SPOp$ loss $\ell_{\SPOp}$ in Section \ref{sec:jpo-non-asymptotic}. We found that the least squares loss satisfies Fisher consistency and uniform calibration under fairly general conditions, but in contrast the $\SPOp$ loss fails to satisfy these conditions in some fairly natural settings. Our numerical results in Section~\ref{sec:num-multiclass} demonstrate that lack of consistency of the loss function can indeed have a detrimental effect on performance. 
}

{\cbl That said, our results in Section \ref{sec:knapsack}  re-affirm \citet{ElmachtoubGrigas2017}'s finding that the $\SPOp$ loss performs well under model misspecification, e.g., when we restrict ourselves to learning a linear predictor but the true underlying data generation model is nonlinear. This suggests a future research direction to build our understanding of robustness to model misspecification of loss functions in the prediction and optimization context. 
Our findings from Sections~\ref{sec:knapsack}~and~\ref{sec:num-multiclass} call for the design of new loss function that are consistent on broad problem classes and take into account optimization problem information as well.
} Some other interesting future directions include further exploration of uniform calibration for loss functions besides  {\cbl $\ell_{\LS}$ and $\ell_{\SPOp}$} and investigating {\cbl theoretical and numerical performance of calibration on} objective functions $f(x,c)$ {\cbl depending non-linearly on $c$}.

\ACKNOWLEDGMENT{This research was supported by NSF grant CMMI 1454548. {\cbl We would like to thank the review team for their suggestions that lead to significant improvements in terms of the presentation of the material.}}

\bibliographystyle{informs2014}

\ECSwitch

\ECHead{Electronic Companion to \emph{Risk Guarantees for End-to-End Prediction and Optimization Processes}}

\section{A Note on the Regularity of $X^*$ and $x^*$}\label{sec:jpo-regularity}

We denote the power set, the collection of all subsets of $X$, as $2^X$. An important property that we exploit is that the argmin mapping $X^*(d)$ is, in a sense, well-behaved as we change $d$. More precisely, the sense of regularity that we use is upper semicontinuity, which stems from a result in perturbation analysis \citep{BonnansShapiro2000}.
\begin{definition}\label{def:jpo-upper-semicty}
A multivalued function $F:\bbR^m \to 2^X$ is \emph{upper semi-continuous} at a point $d \in \bbR^m$ if, for any open set $U$ containing $F(d)$, there exists an open set $U_d$ containing $d$ such that for all $d' \in U_d$, $F(d') \subseteq U$. Equivalently, $F$ is upper semi-continuous if, for any closed set $V$, the following set is closed:
\[ \left\{ d \in \bbR^m : F(d) \cap V \neq \emptyset \right\}. \]
\end{definition}
\begin{lemma}\label{lemma:jpo-X^*-upper-semicty}
Suppose $X$ is compact. Then the multivalued mapping $X^*:\bbR^m \to 2^X$ is upper semi-continuous.
\end{lemma}
\proof{Proof.}
This follows immediately from verifying the conditions of \citet[Proposition 4.4]{BonnansShapiro2000}, which are straightforward to check due to the fact that the domain $X$ does not change with the vector $d$.
\Halmos
\endproof

We can use Lemma \ref{lemma:jpo-X^*-upper-semicty} to show the existence of a measurable selection $x^*(d) \in X^*(d)$ via an application of the Kuratowski--Ryll-Nardzewski theorem on the existence of measurable selectors for multivalued mappings. We use the version stated in \citet[Theorem 6.9.3]{Bogachev2007measure}.
\begin{lemma}\label{lemma:jpo-x^*-measurable}
Suppose $X$ is compact. Then there exists a measurable mapping $x^*:\bbR^m \to X$ such that $x^*(d) \in X^*(d)$ for all $d \in \bbR^m$.
\end{lemma}
\proof{Proof.}
Consider the multivalued function $X^*:\bbR^m \to 2^X$ defined by $X^*(d) = \argmin_{x \in X} d^\top x$. Note that since $d^\top x$ is continuous, $X^*(d) = \{ x \in X : d^\top x = \min_{x' \in X} d^\top x' \}$ is closed (it is the inverse of a singleton). Now consider an open set $U$, and the sets
\[ \hat{X}^*(U) := \left\{ d \in \bbR^m : X^*(d) \cap U \neq \emptyset \right\}. \]
It is known that $U$ can be represented as the countable union of closed sets: $U = \bigcup_{k \in \bbN} V_k$ where $V_k$ are closed. Thus, we can write
\[ \hat{X}^*(U) = \left\{ d \in \bbR^m : \exists k \in \bbN \text{ s.t. } X^*(d) \cap V_k \neq \emptyset \right\} = \bigcup_{k \in \bbN} \left\{ d \in \bbR^m : X^*(d) \cap V_k \neq \emptyset \right\}. \]
Now, since $X^*(d)$ is upper semicontinuous, $\left\{ d \in \bbR^m : X^*(d) \cap U_k \neq \emptyset \right\}$ is closed, hence $\hat{X}^*(U)$ is a countable union of closed sets, hence measurable. This shows that $X^*(\cdot)$ satisfies the conditions of \citet[Theorem 6.9.3]{Bogachev2007measure}, therefore there exists a measurable selection $x^*(d) \in X^*(d)$ for all $d \in \bbR^m$.
\Halmos
\endproof

Furthermore, we can show that \emph{any} selection $x^*$ is at least \emph{Lebesgue} measurable, using the following result of \citet{DrusvyatskiyLewis2011}.
\begin{lemma}[{\citet[Corollary 3.5]{DrusvyatskiyLewis2011}}]\label{lemma:jpo-almost-everywhere-minimizers}
The set
\[ D := \left\{ d \in \bbR^m : X^*(d) \text{ is not a singleton} \right\} \]
has Lebesgue measure zero.
\end{lemma}

\begin{lemma}\label{lemma:jpo-any-x^*-measurable}
Any selection $x^*:\bbR^m \to X$ such that $x^*(d) \in X^*(d)$ for all $d \in \bbR^m$ is Lebesgue measurable.
\end{lemma}
\proof{Proof.}
Lemma \ref{lemma:jpo-x^*-measurable} tells us that there exists one such measurable selection $\bar{x}^*$. Consider another selection $x^*$. Then by Lemma \ref{lemma:jpo-almost-everywhere-minimizers}, $\bar{x}^*$ and $x^*$ differ on at most a set $D$ with Lebesgue measure $0$, which is Lebesgue measurable. Furthermore, all subsets of $D$ are also Lebesgue measurable, so $x^*$ must be Lebesgue measurable.
\Halmos
\endproof

In order for our expectations to be well-defined, we make the following assumption.
\begin{assumption}\label{ass:jpo-lebesgue-measurability}
Any probability distribution $\bbP$ is defined on the $\sigma$-algebra of Lebesgue measurable sets.
\end{assumption}
This is not practically restrictive, since any probability distribution we encounter in practice can be written as a mixture of a distribution which is absolutely continuous with respect to Lebesgue measure (i.e., it has a density function), and a discrete distribution supported on a countable set. Such a probability distribution is Lebesgue measurable.

\section{Proof of Results from Section \ref{sec:jpo-risk}}\label{sec:jpo-proof-risk}

{\ctl
\proof{Proof of Lemma \ref{lemma:jpo-L-nonconvex}.}
Consider two extreme points of $X$, $x_0, x_1$ with $c^\top x_0 > c^\top x_1$. Choose $d_0, d_1$ such that minimizing $d_k^\top x$ over $x \in X$ results in the unique minimum $x_k$ for $k = 0,1$. Now note that $L(d_0,c) - L(d_1,c) = c^\top x_0 - c^\top x_1 > 0$. Let us now consider $d_\gamma =  (1-\gamma) d_0 + \gamma d_1$ for very small $\gamma \in(0,1)$. When $\gamma$ is sufficiently small, then $d_\gamma$ will also have $x_0$ as a unique minimizer, so $L(d_\gamma,c) = L(d_0,c)$. Then because $L(d_0,c) > L(d_1,c)$, we have $L(d_\gamma,c) = L(d_0,c) >(1-\gamma) L(d_0,c) + \gamma L(d_1,c)$. Hence,  $L(d,c)$ is not convex in $d$ for any such $c$.
\Halmos
\endproof
}

\proof{Proof of Lemma \ref{lemma:jpo-true-risk-minimizers}.}
{\cbl The} measurability of $w \mapsto \bbE[c \mid w]$ is obvious by definition of the conditional expectation. Fix some measurable $g:W \to \bbR^m$. Observe that for $w \in W$,
\begin{align*}
	\bbE[L(g(w),c) \mid w] &= \bbE\left[ f(x^*(g(w))) + c^\top x^*(g(w)) - \min_{x \in X} \left\{ f(x) + c^\top x \right\} \mid w \right]\\
	&= {\cbl \bbE\left[f(x^*(g(w))) \mid w \right]} + \bbE\left[ c \mid w \right]^\top x^*(g(w)) - \bbE \left[ \min_{x \in X} \left\{ f(x) + c^\top x \right\} \mid w \right]\\
	&= {\cbl \bbE\left[f(x^*(g(w))) \mid w \right] + g^*(w)^\top x^*(g(w)) - \bbE \left[ \min_{x \in X} \left\{ f(x) + c^\top x \right\} \mid w \right]}\\
	&\geq f(x^*(g^*(w))) + g^*(w)^\top {\cbl x^*(g^*(w))}%
	- \bbE\left[ \min_{x \in X} \left\{ f(x) + c^\top x \right\} \mid w \right]\\
	&= {\cbl f(x^*(g^*(w))) + \bbE\left[ c \mid w \right]^\top x^*(g^*(w))  -  \bbE\left[\min_{x \in X} \left\{ f(x) + c^\top x \right\} \mid w \right]}\\
	&= {\cbl \bbE\left[f(x^*(g^*(w))) + c^\top x^*(g^*(w))  - \min_{x \in X} \left\{ f(x) + c^\top x \right\} \mid w \right]}\\
	&= \bbE[L(g^*(w), c) \mid w],
\end{align*}
where the inequality follows from the definition of $x^*(\cdot)$. 
{\cbl Integrating both sides of this relation over $w \in W$ gives $R(g,\bbP) \geq R(g^*,\bbP)$. Thus, $g^*$ is the minimizer of $R(g,\bbP)$.}

The second result follows because
\begin{align*}
	\min_{d' \in \bbR^m} \bbE[L(d',c) \mid w] &= \min_{d' \in \bbR^m} \left\{f(x^*(d')) + \bbE[c \mid w]^\top x^*(d')\right\} - \bbE\left[ \min_{x \in X} \left\{ f(x) + c^\top x \right\} \mid w \right]\\
	&= {\cbl \min_{d' \in \bbR^m} \left\{f(x^*(d')) + g^*(w)^\top x^*(d')\right\} - \bbE\left[ \min_{x \in X} \left\{ f(x) + c^\top x \right\} \mid w \right]}\\
	&= f(x^*(g^*(w))) + g^*(w)^\top x^*(g^*(w)) - \bbE\left[ \min_{x \in X} \left\{ f(x) + c^\top x \right\} \mid w \right]\\
	&= {\cbl f(x^*(g^*(w))) + \bbE\left[ c \mid w \right]^\top x^*(g^*(w))  -  \bbE\left[\min_{x \in X} \left\{ f(x) + c^\top x \right\} \mid w \right]}\\
	&= {\cbl \bbE\left[ L(g^*(w),c) \mid w \right],}
\end{align*}
{\cbl and then integrating both sides over $w \in W$ gives $\bbE\left[ \min_{d' \in \bbR^m} \bbE[L(d',c) \mid w] \right] = \bbE[L(g^*(w),c)] = R(\bbP)$.}
\Halmos
\endproof

\section{Proof of Theorem \ref{thm:jpo-calibration}}\label{sec:jpo-proof-calibration}

Define
\begin{equation}\label{eqn:jpo-calibration-function}
\delta_\ell(\epsilon,w;\bbP) := \inf_{d \in \bbR^m} \left\{ \bbE[\ell(d,c) \mid w] - \min_{d' \in \bbR^m} \bbE[\ell(d',c) \mid w] : \bbE[L(d,c) \mid w] - \min_{d' \in \bbR^m} \bbE[L(d',c) \mid w] ~{\cbl \geq}~ \epsilon \right\}.
\end{equation}
Note that if $\ell$ is $\bbP$-calibrated, then $\delta_\ell(\epsilon,w;\bbP) > 0$ for all $\epsilon > 0, w \in W$ {\cbl by taking the contrapositive of the implication in Definition \ref{def:jpo-calibration}}. In order to prove Theorem \ref{thm:jpo-calibration}, we first verify measurability for $\delta_{\ell}$.
\begin{lemma}\label{lemma:jpo-delta-measurable}
Suppose $\ell$ is measurable and satisfies Assumption \ref{ass:jpo-finiteness}, and that $X$ is compact. For any $\epsilon > 0$, the function $\delta_\ell(\epsilon, \cdot;\bbP):W \to \bbR$ is measurable.
\end{lemma}
\proof{Proof.}
Consider the set
\[ W_r := \left\{ w \in W : \delta_\ell(\epsilon,w;\bbP) \leq r \right\}. \]
Showing measurability of $\delta_\ell(\epsilon,\cdot;\bbP)$ boils down to showing that $W_r$ is measurable. Rewrite
\begin{align*}
W_r &= \left\{ w \in W : \forall k \in \bbN,\ \exists d \in \bbR^m \text{ s.t. } \begin{aligned}
\bbE[\ell(d,c) \mid w] - \min_{d' \in \bbR^m} \bbE[\ell(d',c) \mid w] &\leq r + 1/k\\
\bbE[L(d,c) \mid w] - \min_{d' \in \bbR^m} \bbE[L(d',c) \mid w] &~{\cbl \geq}~ \epsilon
\end{aligned} \right\}\\
&= \bigcap_{k \in \bbN} \left\{ w \in W : \exists d \in \bbR^m \text{ s.t. } \begin{aligned}
\bbE[\ell(d,c) \mid w] - \min_{d' \in \bbR^m} \bbE[\ell(d',c) \mid w] &\leq r + 1/k\\
\bbE[L(d,c) \mid w] - \min_{d' \in \bbR^m} \bbE[L(d',c) \mid w] &~{\cbl \geq}~ \epsilon
\end{aligned} \right\}
\end{align*}

To this end, first consider the subset
\begin{align*}
W_L(\epsilon) &= \left\{ (w,d) \in W \times \bbR^m : \bbE[L(d,c) \mid w] - \min_{d' \in \bbR^m} \bbE[L(d',c) \mid w] ~{\cbl \geq}~ \epsilon \right\}\\
&= \left\{ (w,d) \in W \times \bbR^m : f(x^*(d)) + \bbE[c \mid w]^\top x^*(d) - \min_{x \in X} \left\{ f(x) + \bbE[c \mid w]^\top x \right\} ~{\cbl \geq}~ \epsilon \right\}.
\end{align*}
This is measurable since $\bbE[c \mid w]$ is measurable in $w$ by definition of conditional expectation, $f$ is continuous hence measurable, and we have assumed $x^*(d)$ is measurable in $d$, which is possible by Lemma \ref{lemma:jpo-x^*-measurable}.

Now consider the subset
\begin{align*}
W_\ell(\alpha) &= \left\{ (w,d) \in W \times \bbR^m : \bbE[\ell(d,c) \mid w] - \min_{d' \in \bbR^m} \bbE[\ell(d',c) \mid w] \leq \alpha \right\}.
\end{align*}
First observe that the function $h$ defined by $h(w,d) = \bbE[\ell(d,c) \mid w]$ is continuous in $d$ and measurable in $w$. Continuity in $d$ follows because $\ell(d,c)$ is convex in $d$, and $h(w,d)$ is finite for any $w$ by Assumption \ref{ass:jpo-finiteness}, and all convex functions are continuous in the relative interiors of their domains (see e.g., \citet[Theorem 10.1]{Rockafellar1970book}). Measurability follows from measurability of $\ell$ and the definition of conditional expectation.

We now show that $h$ is jointly measurable in $(w,d)$ by showing that it is a pointwise limit of measurable functions. For $k \in \bbN$, consider the box $B_k := [-k,k]^m \subset \bbR^m$ and a finite set of grid points $G_k \subset B_k$ such that any point $d \in B_k$ is at most distance $1/k$ away from a grid point in Euclidean norm. If $d \in B_k$, define $h_k(w,d) = h(w,g)$ where $g \in B_k$ is the closest grid point to $d$ (with ties broken arbitrarily), and if $d \not\in B_k$ define $h_k(w,d) = 0$. Note that fixing $g$, $w \mapsto h_g(w) := h(w,g)$ is measurable in $w$. Now, $h_k$ is the sum of finitely many functions of the form $\bm{1}_{D}(d) h_g(w)$ for some measurable set $D$ and grid point $g$. It is easy to check that this is measurable, therefore $h_k$ is measurable. Furthermore, by continuity of $h$ in $d$, $h_k(w,d) \to h(w,d)$ pointwise. Therefore, $h$ is measurable. Finally, the function $(w,d) \mapsto \min_{d' \in \bbR^m} h(w,d')$ is measurable because by continuity of $h$ in $d$, we can write
\[ \left\{ (w,d) : \min_{d' \in \bbR^m} h(w,d') \leq \alpha \right\} = \bigcup_{d' \in D_*, k \in \bbN} \left\{ (w,d) : h(w,d') \leq \alpha + 1/k \right\} \]
where $D_*$ is a countable dense subset of $\bbR^m$ (e.g., $\bbQ^m$). This shows that $W_\ell(\alpha)$ is measurable because the function $h(w,d) - \min_{d' \in \bbR^m} h(w,d') = \bbE[\ell(d,c) \mid w] - \min_{d' \in \bbR^m} \bbE[\ell(d',c) \mid w]$ is measurable.

Now notice that the set
\[ \left\{ (w,d) \in W \times \bbR^m : \begin{aligned}
\bbE[\ell(d,c) \mid w] - \min_{d' \in \bbR^m} \bbE[\ell(d',c) \mid w] &\leq r + 1/k\\
\bbE[L(d,c) \mid w] - \min_{d' \in \bbR^m} \bbE[L(d',c) \mid w] &~{\cbl \geq}~ \epsilon
\end{aligned} \right\} = W_\ell(r+1/k) \cap W_L(\epsilon) \]
is measurable. Therefore, its projection onto $W$ is measurable, which is
\[\left\{ w \in W : \exists d \in \bbR^m \text{ s.t. } \begin{aligned}
\bbE[\ell(d,c) \mid w] - \min_{d' \in \bbR^m} \bbE[\ell(d',c) \mid w] &\leq r + 1/k\\
\bbE[L(d,c) \mid w] - \min_{d' \in \bbR^m} \bbE[L(d',c) \mid w] &~{\cbl \geq}~ \epsilon
\end{aligned} \right\}.\]
This shows that $W_r$ is measurable, concluding our proof.
\Halmos
\endproof

\proof{Proof of Theorem \ref{thm:jpo-calibration}.}
{\cbl We wish to apply a result of \citet[Theorem 2.8]{Steinwart2007}, for which we need to show that there exists measurable functions $b:W \to \bbR$ and and $\delta:(0,\infty)\times W \to (0,\infty)$ such that $\bbE[|b(w)|] < \infty$, for any $w \in W$
\[\bbE[L(d,c) \mid w] - \min_{d' \in \bbR^m} \bbE[L(d',c) \mid w] \leq b(w),\]
and for any $\epsilon > 0$ and predictor $g:W \to \bbR^m$,
\begin{align*}
&\left\{ w \in W : \bbE[\ell(g(w),c) \mid w] - \min_{d' \in \bbR^m} \bbE[\ell(d',c) \mid w] < \delta(\epsilon,w) \right\}\\
&\subseteq \left\{ w \in W : \bbE[L(g(w),c) \mid w] - \min_{d' \in \bbR^m} \bbE[L(d',c) \mid w] < \epsilon \right\}.
\end{align*}

We first find $b$. Let $\Omega$ be the $\ell_\infty$-diameter of the set $X$, which is finite since $X$ is compact. Observe that for any $d \in \bbR^m$, 
\begin{align*}
	\bbE[L(d,c) \mid w] - \min_{d' \in \bbR^m} \bbE[L(d',c) \mid w] &= f(x^*(d)) + \bbE[c \mid w]^\top x^*(d) - \min_{x \in X} \left\{ f(x) + \bbE[c \mid w]^\top x\right\}\\
	&\leq \max_{x,x' \in X} \left\{ f(x) - f(x') + \bbE[c \mid w]^\top (x'-x) \right\}\\
	&\leq \max_{x,x' \in X} \left\{ f(x) - f(x') + \|\bbE[c \mid w]\|_1 \|x'-x\|_\infty \right\}\\
	&\leq \Omega \|\bbE[c \mid w]\|_1 + \max_{x,x' \in X} \left\{ f(x) - f(x') \right\}.
\end{align*}
Therefore, we can define $b(w) := \Omega \|\bbE[c \mid w]\|_1 + \max_{x,x' \in X} \left\{ f(x) - f(x') \right\}$ for each $w \in W$, which is integrable as $\bbE\left[ \|\bbE[c \mid w]\|_1 \right] \leq \bbE\left[ \bbE[\|c\|_1 \mid w] \right] = \bbE[\|c\|_1] < \infty$ by Assumption \ref{ass:jpo-finiteness}.

We will take $\delta := \delta_\ell(\cdot; \bbP)$ defined in \eqref{eqn:jpo-calibration-function}, which is measurable by Lemma \ref{lemma:jpo-delta-measurable}. For any $g:W \to \bbR^m$, and $w \in W$ such that $\bbE[L(g(w),c) \mid w] - \min_{d' \in \bbR^m} \bbE[L(d',c) \mid w] \geq \epsilon$, by $\bbP$-calibration and definition of $\delta_{\ell}$ we have $\bbE[\ell(g(w),c) \mid w] - \min_{d' \in \bbR^m} \bbE[\ell(d',c) \mid w] \geq \delta_\ell(\epsilon,w;\bbP)$, therefore the required property for $\delta$ is satisfied.

Applying the result of \citet[Theorem 2.8]{Steinwart2007} then gives the risk bound.}
\Halmos
\endproof

\section{Proofs of Results from Section \ref{sec:jpo-FisherConsistent-surrogate-risk}}\label{sec:jpo-proofs-FisherConsistent}

\proof{Proof of Theorem \ref{thm:jpo-FisherConsistent-calibration}.}
Denote
\begin{align*}
D_\ell(\alpha;w) &:= \left\{ d \in \bbR^m : \bbE[\ell(d,c) \mid w] - \min_{d' \in \bbR^m} \bbE[\ell(d,c) \mid w] ~{\cbl <}~ \alpha \right\}\\
D(\alpha;w) &:= \left\{ d \in \bbR^m : \bbE[L(d,c) \mid w] - \min_{d' \in \bbR^m} \bbE[L(d,c) \mid w] ~{\cbl <}~ \alpha \right\}.
\end{align*}
Note that
\[ \argmin_{d' \in \bbR^m} \bbE[\ell(d',c) \mid w] = \bigcap_{\alpha > 0} D_\ell(\alpha;w), \quad \argmin_{d' \in \bbR^m} \bbE[L(d',c) \mid w] = \bigcap_{\alpha > 0} D(\alpha;w). \]
Suppose first that $\ell$ is $\bbP$-calibrated. Then for any $\epsilon > 0$, there exists $\delta > 0$ (which can depend on $w$) such that $D_\ell(\delta;w) \subseteq D(\epsilon;w)$. In particular, since $D_\ell(\alpha;w) \subseteq D_\ell(\alpha';w)$ for $\alpha \leq \alpha'$, we have
\[ \argmin_{d' \in \bbR^m} \bbE[\ell(d',c) \mid w] = \bigcap_{0 < \alpha \leq \delta} D_\ell(\alpha;w) \subseteq D(\epsilon;w). \]
Taking the intersection of the right hand side over $\epsilon > 0$, we have
\[ \argmin_{d' \in \bbR^m} \bbE[\ell(d',c) \mid w] \subseteq \argmin_{d' \in \bbR^m} \bbE[L(d',c) \mid w], \]
hence $\ell$ is $\bbP$-Fisher consistent.

Suppose now that $\ell$ is not $\bbP$-calibrated. We show that it is also not $\bbP$-Fisher consistent. Fix an arbitrary $w \in W$. Note that the function $h:\bbR^m \to \bbR$ defined by $h(d) = \bbE[\ell(d,c) \mid w]$ is convex by convexity of $\ell(d,c)$, and hence under Assumption \ref{ass:jpo-finiteness}, it is continuous (see e.g., \citet[Theorem 10.1]{Rockafellar1970book}). 

Since $\ell$ is not $\bbP$-calibrated, there exists $w \in W$ and $\epsilon > 0$ such that for all $\delta > 0$, there exists $d(\delta) \in \bbR^m$ such that $h(d(\delta)) - \min_{d' \in \bbR^m} h(d') ~{\cbl <}~ \delta$ but $\bbE[L(d(\delta),c) \mid w] - \min_{d' \in \bbR^m} \bbE[L(d,c) \mid w] ~{\cbl \geq}~ \epsilon$.

Now, let $d_k = d(1/k)$ for $k \in \bbN$. Note that $\{d_k\}_{k \in \bbN} \subset D_\ell(1;w)$ which is compact since by Assumption \ref{ass:jpo-finiteness} $\argmin_{d' \in \bbR^m} h(d')$ is compact, so all level sets are bounded (see, e.g., \citet[Corollary 8.7.1]{Rockafellar1970book}). Therefore, there exists a convergent subsequence $d_k' \to d \in {\cbl \cl D_\ell(1;w)}$. Since $h$ is continuous, we must have $d \in \argmin_{d' \in \bbR^m} h(d')$.

We now want to show that $d \not\in \argmin_{d' \in \bbR^m} \bbE[L(d',c) \mid w]$. We know from Lemma \ref{lemma:jpo-X^*-upper-semicty} that the argmin mapping $X^*(\cdot)$ is upper semi-continuous at $d$. Suppose for contradiction that $d \in \argmin_{d' \in \bbR^m} \bbE[L(d',c) \mid w]$. Then we must have $X^*(d) \subseteq X^*(\bbE[c \mid w])$. Thus, for $\epsilon > 0$ the set
\[ X^\circ(\epsilon') = \left\{ x' : f(x') + \bbE[c \mid w]^\top x' < \min_{x \in X} \left\{ f(x) + \bbE[c \mid w]^\top x \right\} + \epsilon' \right\} \]
is an `open' set (as $x \mapsto f(x) + \bbE[c \mid w]^\top x$ is continuous) containing $X^*(d)$. Note that this is not open in $\bbR^m$ by the usual topology, since $f(x)$ may be infinite for $x \not\in X$. However, it is open when we work with $X \subset \bbR^m$ as the \emph{entire} topological space with the induced topology from $\bbR^m$. Then, by Definition \ref{def:jpo-upper-semicty} of upper semi-continuity, there exists a neighbourhood $D^\circ(\epsilon')$ of $d$ such that for any $d^\circ \in D^\circ(\epsilon')$, $X^*(d^\circ) \subset X^\circ(\epsilon')$, which means that $\bbE[L(d^\circ,c) \mid w] < \min_{d' \in \bbR^m} \bbE[L(d',c) \mid w] + \epsilon'$ since $x^*(d^\circ) \in X^*(d^\circ) \subseteq X^\circ(\epsilon')$.

But now consider $\epsilon' < \epsilon$. Since $d_k' \to d$, $D^\circ(\epsilon')$ is open, and $d \in D^\circ(\epsilon')$, we eventually have $d_k' \in D^\circ(\epsilon')$ for sufficiently large $k$. But this contradicts the fact that by construction of the sequence $\{d_k\}_{k  \in \bbN}$ we have $\epsilon' < \epsilon < \bbE[L(d_k',c) \mid w] - \min_{d' \in \bbR^m} \bbE[L(d',c) \mid w] = \bbE[c \mid w]^\top x^*(d_k') - \min_{x \in X} \bbE[c \mid w]^\top x$.
\Halmos
\endproof

\proof{Proof of Corollary \ref{cor:jpo-FisherConsistent-conv}.}
Fix some $\epsilon > 0$. Take $\delta > 0$ corresponding to $\epsilon$ in Corollary \ref{cor:jpo-FisherConsistent}. Since $R_\ell(g_n,\bbP) \to R_\ell(\bbP)$, we have $R_\ell(g_n,\bbP) \leq R_\ell(\bbP) + \delta$ eventually. By Theorem \ref{cor:jpo-FisherConsistent}, we will also have $R(g_n,\bbP) \to R(\bbP) + \epsilon$ eventually. 
\Halmos
\endproof

\proof{{\cbl Proof of Example \ref{ex:jpo-spo-one-dim}}.}
Let us explore what $\argmin_{d' \in \bbR} \bbE[\ell(d,c) \mid w]$ is for our setting. For convenience, we fix $w \in W$, and omit the $w$ in the notation, so that $D^* = D_w^*$, $\bbE[\cdot] = \bbE[\cdot \mid w]$ and $\bbP[\cdot] = \bbP[\cdot \mid w]$. Then
\[ 2\bbE[\ell(d,c)] = \bbE[|2d-c|] - 2d \bbE[\sign(c)] + \bbE[|c|] = \bbE[|2d-c|] + 2d \left(\bbP[c < 0] - \bbP[c > 0]\right) + \bbE[|c|]. \]
This is a convex function in $d$, so we look at the subdifferential to determine its minimizers. Note that
\[ \partial_d \bbE[|2d-c|] = \left\{ 2\left( \bbP[c < 2d] - \bbP[c > 2d] + s \bbP[c=2d] \right) : s \in [-1,1] \right\}, \]
so
\[ \partial_d \bbE[\ell(d,c)] = \left\{ \bbP[c < 2d] - \bbP[c > 2d] + \bbP[c < 0] - \bbP[c > 0] + s \bbP[c=2d] : s \in [-1,1] \right\}. \]
For simplicity, let us assume that $\bbP[c=2d] = 0$ for any $d$ (many such distributions exist). Then $\bbE[\ell(d,c)]$ is differentiable with
\[ \grad_d \bbE[\ell(d,c)] = \bbP[c < 2d] - \bbP[c > 2d] + \bbP[c < 0] - \bbP[c > 0]. \]
Denote $d^*$ to be a minimizer of $\bbE[\ell(d,c)]$. If $\bbP[c < 0] = \bbP[c > 0]$, then setting $d=0$ gives $\grad_d \bbE[\ell(d,c)] = 0$, so $d^* = c = 0$. If $\bbP[c < 0] - \bbP[c > 0] < 0$, then $\left. \grad_d \bbE[\ell(d,c)] \right|_{d=0} < 0$, so increasing $d$ from $0$ will decrease $\bbE[\ell(d,c)]$. Thus, $d^* > 0$. However, note that $\bbP[c < 0] - \bbP[c > 0] < 0$ implies that the median of $c$ is also $>0$. If $\bbP[c < 0] - \bbP[c > 0] > 0$, then $\left. \grad_d \bbE[\ell(d,c)] \right|_{d=0} > 0$, so decreasing $d$ from $0$ will decrease $\bbE[\ell(d,c)]$. Thus, $d^* < 0$. However, note that $\bbP[c < 0] - \bbP[c > 0] > 0$ implies that the median of $c$ is also $<0$. In all cases, the minimizer $d^*$ is of the same sign as the median of $c$. Now, if $\bbP$ is a symmetric distribution, then the mean $\bbE[c]$ is equal to the median, and thus $d^*$ has the same sign as $\bbE[c]$, so also minimizes $\bbE[L(d,c)]$. However, if the median has a different sign to the mean, then $\ell$ is not $\bbP$-Fisher consistent. Such distributions can be constructed by shifting a log-normal distribution, for example.
\Halmos
\endproof

\proof{Proof of Example \ref{ex:jpo-multiclass}.}
With the distribution $\bbP$ specified, $\min_{d' \in \bbR^m} \bbE[\ell_{\SPOp}(d',c)]$ can be expressed as the following linear program (making the change of variables $2d' \to d$):
\begin{align*}
\min_{d,\gamma} &\quad \sum_{j \in [m]} p_j (d_j - \gamma_j)\\
\text{s.t.} &\quad \gamma_j \leq d_j, \ j \in [m]\\
&\quad \gamma_j \leq d_k-1, \ j,k \in [m], k \neq j\\
&\quad d,\gamma \in \bbR^m.
\end{align*}

We analyse this linear program. Fix a vector $d \in \bbR^m$. Let $j^* \in \argmin_{j' \in [m]} d_{j'}$. Then since $p_k > 0$ for all $k \neq j^*$, the optimal choice of $\gamma_k$ makes it as large as possible, so we set $\gamma_k = d_{j^*}-1$ for $k \neq j^*$. In other words, for all but one index $j^* \in \argmin_{j' \in [m]} d_{j'}$, we set $\gamma_j = \min_{j' \in [m]} d_{j'} - 1$. For $j^*$, we set $\gamma_{j^*} = \min\left\{ d_{j^*}, \min_{j' \neq j^*} d_{j'} - 1 \right\}$.

If there exists $j \neq j^*$ such that $d_{j^*} \leq d_j - 1$, then decreasing $d_j \downarrow d_{j^*}+1$ does not violate any constraints since $\gamma_j = d_{j^*}-1 < d_j$ and $\gamma_{j^*} = d_{j^*} \leq d_j - 1$, and decreases the objective. Therefore, without loss of generality, we assume that $d_j-1 \leq d_{j^*}$ for all $j \neq j^*$. This implies that $\gamma_{j^*} = \min_{j' \neq j^*} d_{j'} - 1$.

Furthermore, if we have $j,k \in [m] \setminus \{j^*\}$, $j \neq k$ such that $d_j < d_k$, note that we can decrease $d_k \downarrow d_j$ without violating any constraints, since $\gamma_{j'} = d_{j^*} - 1 \leq d_j - 1 < d_k - 1 < d_k$ for all $j' \neq j^*$ and $\gamma_{j^*} \leq d_j - 1 < d_k - 1$. This implies that, without loss of generality, we can assume that for $j \neq j^*$, we have $d_j = \delta$ for some $\delta \in [d_{j^*}, d_{j^*} + 1]$. In particular, this implies that $\gamma_{j^*} = \delta - 1$, thus the objective becomes
\[ \sum_{j \in [m]} p_j (d_j - \gamma_j) = (\delta - d_{j^*} + 1) \sum_{j \neq j^*} p_j + p_{j^*} (d_{j^*} - \delta + 1) = (1-2p_{j^*}) (\delta - d_{j^*}) + 1. \]
This shows that if $p_{j^*} > 1/2$, then we should make $\delta$ as large as possible, i.e., $\delta = d_{j^*} + 1$. On the other hand, when $p_{j^*} < 1/2$, we set $\delta = d_{j^*}$, i.e., the optimal vector $d^*$ is constant.

This implies that, if there exists $j^* \in [m]$ such that $p_{j^*} > 1/2$, and necessarily $j^* = \argmax_{j' \in [m]} p_{j'}$, then the minimizers of $\bbE[\ell(d,c)]$ take the form $d_\alpha = (\alpha \bm{1}_m - e_{j^*})/2$ for $\alpha \in \bbR$. Clearly, $\argmin_{j' \in [m]} d_{\alpha,j'} = j^*$, so for such distributions $\bbP$, $\ell_{\SPOp}$ is $\bbP$-Fisher consistent.

On the other hand, for distributions $\bbP$ with $\max_{j' \in [m]} p_{j'} < 1/2$, $\ell_{\SPOp}$ is not $\bbP$-Fisher consistent, since the set of minimizers of $\bbE[\ell(d,c)]$ are the vectors $d_\alpha = \alpha \bm{1}_m$, $\alpha \in \bbR$, which cannot in general pick out the maximum probability class $j \in [m]$, i.e., the highest $p_j$.
\Halmos
\endproof

\section{Proofs of Results from Section \ref{sec:jpo-non-asymptotic}}\label{sec:jpo-proofs-non-asymptotic}

\proof{Proof of Lemma \ref{lemma:jpo-uniform-calibration}.}
{\cbl The `only if' direction was established in Remark \ref{rem:jpo-uniform-calibration-function}, so we only need to prove the `if' direction.}

When $\delta_\ell(\epsilon; \CP) > 0$, take $0 < \delta \leq \delta_\ell(\epsilon; \CP)$, and noting that $\delta_\ell(\cdot; \CP)$ is non-decreasing, we get for any $d \in \bbR^m$, $w \in W$ and $\bbP \in \CP$,
\[\bbE[\ell(d,c) \mid w] - \min_{d' \in \bbR^m} \bbE[\ell(d',c) \mid w] \leq \delta < \delta_\ell(\epsilon; \CP).\]
If $d \in \bbR^m$, $w \in W$ and $\bbP \in \CP$ were such that $\bbE[L(d,c) \mid w] - \min_{d' \in \bbR^m} \bbE[L(d',c) \mid w] > \epsilon$, we reach a contradiction since we would then by definition of $\delta_\ell(\cdot; \CP)$ in \eqref{eqn:jpo-uniform-calibration-function} have $\bbE[\ell(d,c) \mid w] - \min_{d' \in \bbR^m} \bbE[\ell(d',c) \mid w] \geq \delta_\ell(\epsilon;\CP)$. Thus, for any $w \in W$ and $\bbP \in \CP$,
\[\bbE[\ell(d,c) \mid w] - \min_{d' \in \bbR^m} \bbE[\ell(d',c) \mid w] \leq \delta \implies \bbE[L(d,c) \mid w] - \min_{d' \in \bbR^m} \bbE[L(d',c) \mid w] \leq \epsilon. \]
\Halmos
\endproof

\proof{Proof of Theorem \ref{thm:jpo-uniform-calibration}.}
{\cbl When $\ell$ is $\cP$-uniformly calibrated, we know that $\delta_{\ell}(\epsilon;\cP) > 0$ for any $\epsilon > 0$. \citet[Lemma A.6]{Steinwart2007} shows that this implies $\delta^{**}(\epsilon;\cP) > 0$ for $\epsilon \in (0,B_f + B_C B_X]$.

We can now utilize \citet[Theorem 2.13]{Steinwart2007} to derive the risk bound. In order to do so, note that \citet[Theorem 2.13]{Steinwart2007} requires us to verify that for any $g:W \to \bbR$,
\[ \esssup_{w \in W} \bbE[L(g(w),c) \mid w] - \min_{d' \in \bbR^m} \bbE[L(d',c) \mid w] \leq B_f + B_C B_X, \]
where $\esssup$ stands for essential supremum. 
The relation above follows from Remark \ref{rem:jpo-boundedness} and from the definition of $\delta_{\ell}(\cdot;\cP)$ that ensures that for any $\epsilon > 0$ and $w \in W$, we have
\[ \left\{ g : \bbE[\ell(g(w),c) \mid w] - \min_{d' \in \bbR^m} \bbE[\ell(d',c) \mid w] < \delta_{\ell}(\epsilon; \cP) \right\} \subseteq \left\{ g : \bbE[L(g(w),c) \mid w] - \min_{d' \in \bbR^m} \bbE[L(d',c) \mid w] < \epsilon \right\}. \]
}
\Halmos
\endproof

\proof{Proof of Lemma \ref{lemma:jpo-uniform-calibration-function}.}
Fix an arbitrary $w \in W$. Note that $\bbE[c \mid w] \in \{c' : x^*(c') \in X^*(\bbE[c \mid w])\} = \argmin_{d' \in \bbR^m} \bbE[L(d',c) \mid w]$, hence we have
\begin{align*}
\bbE[L(d,c) \mid w] - \min_{d' \in \bbR^m} \bbE[L(d',c) \mid w] &= f(x^*(d)) + \bbE[c \mid w]^\top x^*(d) - \min_{d' \in \bbR^m} \left\{ f(x^*(d')) + \bbE[c \mid w]^\top x^*(d') \right\}\\
&= f(x^*(d)) - f(x^*(\bbE[c \mid w])) + \bbE[c \mid w]^\top (x^*(d) - x^*(\bbE[c \mid w])).
\end{align*}
Hence
\begin{align*}
\delta_\ell(\epsilon;\CP) &= \inf\limits_{\substack{d \in \bbR^m\\ w \in W\\ \bbP \in \CP}} \left\{ \bbE[\ell(d,c) \mid w] - \min_{d' \in \bbR^m} \bbE[\ell(d',c) \mid w] : \bbE[L(d,c) \mid w] - \min_{d' \in \bbR^m} \bbE[L(d',c) \mid w] ~{\cbl \geq}~ \epsilon \right\}\\
&= \inf\limits_{\substack{d,\bar{c} \in \bbR^m}} \inf\limits_{\substack{w \in W \\ \bbP \in \CP \\ \bbE[c \mid w] = \bar{c}}} \left\{ \bbE[\ell(d,c) \mid w] - \min_{d' \in \bbR^m} \bbE[\ell(d',c) \mid w] : f(x^*(d)) - f(x^*(\bar{c})) + \bar{c}^\top (x^*(d) - x^*(\bar{c})) ~{\cbl \geq}~ \epsilon \right\}\\
&= \inf\limits_{x,x' \in X} \inf\limits_{\substack{d : x^*(d) = x\\ \bar{c} : x^*(\bar{c}) = x'}} \inf\limits_{\bbP : \bbE[c] = \bar{c}} \left\{ \bbE[\ell(d,c)] - \min_{d' \in \bbR^m} \bbE[\ell(d',c)] : f(x) - f(x') + \bar{c}^\top (x - x') ~{\cbl \geq}~ \epsilon \right\}.
\end{align*}
\Halmos
\endproof

\proof{Proof of Lemma \ref{lemma:jpo-uniform-calibration-distance-bound}.}
Fix arbitrary distinct $x,x' \in X$. Consider the halfspace
\begin{align*}
H_0(x,x') = \left\{ d' : f(x) - f(x') + (d')^\top (x - x') \leq 0 \right\} \supseteq \{d' : x^*(d') = x \}.
\end{align*}
Since $x^*(d) = x$, we have $d \in H_0(x,x')$. Now, if $\bar{c} \in H_0(x,x')$, then $f(x) - f(x') + \bar{c}^\top (x - x') \leq 0$, hence $\|d - \bar{c}\|_2 \geq 0$.

On the other hand, if $x^*(\bar{c}) = x'$ and $f(x) - f(x') + \bar{c}^\top (x - x') > 0$, then we have $\bar{c} \not\in H_0(x,x')$, hence the distance between $\bar{c}$ and $d$ is bounded below by the distance between $\bar{c}$ and the halfspace $H_0(x,x')$, which has the expression
\[ \|d - \bar{c}\|_2 \geq \inf_{d' \in H_0(x,x')} \|d' - \bar{c}\|_2 = \frac{f(x) - f(x') + \bar{c}^\top (x - x')}{\|x - x'\|_2}. \]
\Halmos
\endproof

\proof{Proof of Lemma \ref{lemma:jpo-delta-simple-squared-loss}.}
The usual bias-variance decomposition for squared error gives us
\begin{align*}
\bbE[\ell_{\LS}(d,c) \mid w] &= \bbE[\|d-c\|_2^2 \mid w]\\
&= \|d - \bbE[c \mid w]\|_2^2 + 2\bbE\left[ (d-\bbE[c \mid w])^\top (\bbE[c \mid w] - c)\right] + \bbE\left[ \left\| \bbE[c \mid w] - c \right\|_2^2 \right]\\
&= \|d - \bbE[c \mid w]\|_2^2 + \bbE\left[ \left\| \bbE[c \mid w] - c \right\|_2^2 \right].
\end{align*}
Hence, we can minimize this by choosing $d = \bbE[c \mid w]$, and
\[ \min_{d' \in \bbR^m} \bbE[\ell_{\LS}(d',c) \mid w] = \bbE\left[ \left\| \bbE[c \mid w] - c \right\|_2^2 \right]. \]
Therefore,
\[ \bbE[\ell_{\LS}(d,c) \mid w] - \min_{d' \in \bbR^m} \bbE[\ell_{\LS}(d',c) \mid w] = \|d - \bbE[c \mid w]\|_2^2. \]
Substituting this into \eqref{eqn:jpo-uniform-calibration-alternate} and using the fact that the definition of $\CP$ tells us that $\bar{c} = \bbE[c \mid w]$ can take on any point in $\Conv(C)$ gives the result.
\Halmos
\endproof

\proof{Proof of Theorem \ref{thm:jpo-delta-squared-loss}.}
Fixing distinct $x,x' \in X$, notice that if $\bar{c}, d$ are chosen according to the conditions of Lemma \ref{lemma:jpo-uniform-calibration-distance-bound}, together with the condition that $f(x) - f(x') + \bar{c}^\top (x - x') > \epsilon$, then $\|d - \bar{c}\|_2 > \epsilon/\|x - x'\|_2 \geq \epsilon/B_X > 0$. Together with Lemma \ref{lemma:jpo-delta-simple-squared-loss}, we have, for all $\epsilon > 0$,
\[ \delta_{\ell_{\LS}}(\epsilon;\CP) \geq \frac{\epsilon^2}{B_X^2} > 0. \]
Then $\CP$-uniform calibration follows from Lemma \ref{lemma:jpo-uniform-calibration}.
\Halmos
\endproof

\proof{Proof of Corollary \ref{cor:jpo-squared-loss-risk-bound}.}
The result follows by observing that $\epsilon^2/B_X^2 \leq \delta^{**}(\epsilon)$ since $\epsilon \mapsto \epsilon^2/B_X^2$ is already convex, and then applying Theorem \ref{thm:jpo-uniform-calibration}.
\Halmos
\endproof

\proof{Proof of Theorem \ref{thm:jpo-squared-loss-calibrate}.}
Analogous to \eqref{eqn:jpo-uniform-calibration-function}, define
\[ \delta_j(\epsilon;\CP_{\sym}) := \inf\limits_{\substack{d_j \in \bbR\\ w \in W\\ \bbP \in \CP_{\sym}}} \left\{ \bbE[\ell_j(d_j,c_j) \mid w] - \min\limits_{d_j' \in \bbR} \bbE[\ell_j(d_j',c_j) \mid w] : (d_j - \bbE[c_j \mid w])^2 > \epsilon \right\}. \]
We first show that $\delta_j(\epsilon;\CP_{\sym}) > 0$ for all $\epsilon > 0$.

First, fix $\bbP \in \CP_{\sym}$ and $w \in W$, and observe that for any $d_j$
\begin{align*}
\bbE[\ell_j(\bbE[c_j \mid w] + d_j,c_j) \mid w] &= \bbE[\psi_j(d_j - (c_j - \bbE[c_j \mid w]) \mid w] = \bbE[\psi_j(d_j - (\bbE[c_j \mid w] - c_j)) \mid w]\\
&= \bbE[\psi_j(-d_j + \bbE[c_j \mid w] - c_j) \mid w]\\
&= \bbE[\ell_j(\bbE[c_j \mid w] - d_j,c_j) \mid w].
\end{align*}
Since $\psi$ is strictly convex, $\bbE[\ell_j(d_j,c_j) \mid w]$ is strictly convex in $d_j$, thus for any $d_j \neq 0$,
\begin{align*}
&\bbE[\ell_j(\bbE[c_j \mid w] + d_j,c_j) \mid w] - \bbE[\ell_j(\bbE[c_j \mid w], c_j) \mid w]\\
&= \frac{1}{2} \bbE[\ell_j(\bbE[c_j \mid w] + d_j,c_j) \mid w] + \frac{1}{2} \bbE[\ell_j(\bbE[c_j \mid w] - d_j,c_j) \mid w] - \bbE[\ell_j(\bbE[c_j \mid w], c_j) \mid w]\\
&\geq \frac{1}{4} \bbE\left[ \delta_j(4 d_j^2) \mid w \right] = \frac{1}{4} \delta_j(4d_j^2).
\end{align*}
This shows that $\delta_j(\epsilon;\CP_{\sym}) \geq \delta_j(2\epsilon)/4$. Now, following the outline in Section \ref{sec:jpo-non-asymptotic-outline} and proceeding similarly to the proof of Theorem \ref{thm:jpo-uniform-calibration}, we deduce the risk bound.
\Halmos
\endproof

\proof{Proof of Proposition \ref{prop:jpo-spo-one-dim-non-calibrated}.}
The proof is by construction. We will fix the mean of our class to be $\epsilon$. Let $\phi$ be the density function of the standard normal distribution, and $\Phi$ be the distribution function (note that $\phi(c-\epsilon)$ is the density function of a $N(\epsilon,1)$ random variable. Let $z_\epsilon = \Phi(-\epsilon) + 1 - \Phi(\epsilon)$ denote the probability that a standard normal variable is $<-\epsilon$ or $> \epsilon$. Furthermore, let $\{h(\cdot; \alpha)\}_{\alpha \in (0,1)}$ be a class of continuous functions such that for each $\alpha \in (0,2/3)$, $h(r; \alpha) > 0$ for $r \in [0,1]$, $h(1; \alpha) = 1$, $\int_{r=0}^1 h(r; \alpha) dr = \alpha$. Such a class can be defined as follows:
\[ h(r;\alpha) = \begin{cases}
\alpha/2, & 0 \leq r \leq (2-3\alpha)/(2-\alpha)\\
(2-\alpha)^2(r-1)/(4\alpha) + 1, &(2-3\alpha)/(2-\alpha) < r \leq 1.
\end{cases} \]

For each $k \in \bbN$, define the following density function $\psi^{(k)}$:
\[ \psi^{(k)}(c) = \left\{ \begin{aligned}
&\frac{z_\epsilon + (1-1/k)(1-z_\epsilon)}{z_\epsilon} \phi(c - \epsilon), &c \leq 0\\
&\frac{z_\epsilon + (1-1/k)(1-z_\epsilon)}{z_\epsilon} \phi(-\epsilon) h\left(1-c/\epsilon; \frac{(1-z_\epsilon)z_\epsilon}{2k \epsilon \left( z_\epsilon + (1-1/k) (1-z_\epsilon) \right) \phi(-\epsilon)} \right), &0 < c < \epsilon\\
&\frac{z_\epsilon + (1-1/k)(1-z_\epsilon)}{z_\epsilon} \phi(-\epsilon) h\left(c/\epsilon - 1; \frac{(1-z_\epsilon)z_\epsilon}{2k \epsilon \left( z_\epsilon + (1-1/k) (1-z_\epsilon) \right) \phi(-\epsilon)} \right), &\epsilon \leq c < 2\epsilon\\
&\frac{z_\epsilon + (1-1/k)(1-z_\epsilon)}{z_\epsilon} \phi(\epsilon - c), &c \geq 2\epsilon.
\end{aligned} \right. \]
By construction, $\psi^{(k)}(c)$ is continuous and positive for all $c \in \bbR$, and integrates to $1$. Let $\bbP^{(k)}$ denote the corresponding probability distribution, and by construction we have $\bbP^{(k)}[0 \leq c \leq 2\epsilon] = (1-z_\epsilon)/k$. Therefore $\bbP^{(k)}[c > 0] - \bbP^{(k)}[c < 0] = (1-z_\epsilon)/k \to 0$ as $k \to \infty$, but $\bbE^{(k)}[c] = \epsilon$ since $\psi^{(k)}(c - \epsilon) = \psi(\epsilon - c)$ is symmetric about $\epsilon$.
\Halmos
\endproof

\proof{Proof of Proposition \ref{prop:jpo-spo-one-dim-calibrated}.}
We know that
\[ \delta_{\ell_{\SPOp}}(\epsilon;\CP_{\cont,\sym,\alpha}) = \inf_{w \in W} \inf_{\substack{\bbP \in \CP_{\cont,\sym,\alpha}\\ |\bbE[c \mid w]| > \epsilon}} \left\{ \bbE[c \mid w] \left( \bbP[c > 0 \mid w] - \bbP[c < 0 \mid w] \right) \right\}. \]
Using the property of $\CP_{\cont,\sym,\alpha}$, we deduce that
\[ \delta_{\ell_{\SPOp}}(\epsilon;\CP_{\cont,\sym,\alpha}) \geq \inf_{w \in W} \inf_{\substack{\bbP \in \CP_{\cont,\sym,\alpha}\\ |\bbE[c \mid w]| > \epsilon}} \alpha |\bbE[c \mid w]| \geq \alpha \epsilon. \]
Thus since $\delta_{\ell_{\SPOp}}(\epsilon;\CP_{\cont,\sym,\alpha}) > 0$ for any $\epsilon > 0$, we have uniform calibration by Lemma \ref{lemma:jpo-uniform-calibration}. Furthermore, Theorem \ref{thm:jpo-uniform-calibration} gives us the risk bound.
\Halmos
\endproof

{\cbl
\section{Proof of Results from Section \ref{sec:jpo-computational}}\label{sec:jpo-proofs-computational}

\proof{Proof of Proposition \ref{prop:jpo-portfolio-SPOp-longshort}.}
Using Lagrange duality we know that the dual problem is
\[ \min_{\gamma} \left\{ b \gamma + \frac{1}{2} (d-\gamma p)^\top Q^{-1} (d-\gamma p) \right\} = -\min_x \left\{ \frac{1}{2} x^\top Q x - d^\top x : p^\top x = b \right\}, \]
and the optimal solution is $x^*(d) = Q^{-1} (d - \gamma^* p)$ where $\gamma^*$ is the optimal dual solution. The closed form solution is $\gamma^* = \frac{1}{p^\top Q^{-1} p} \left(p^\top Q^{-1} d - b \right)$, hence
\[ x^*(d) = Q^{-1} d - \frac{1}{p^\top Q^{-1} p} \left(p^\top Q^{-1} d - b \right) Q^{-1} p = A d + \frac{b}{p^\top Q^{-1} p} Q^{-1} p. \]
Observe that $Ap = 0$, so
\begin{align*}
	x^*(d)^\top Q x^*(d) &= \left( A d + \frac{b}{p^\top Q^{-1} p} Q^{-1} p \right)^\top \left( Q A d + \frac{b}{p^\top Q^{-1} p} p \right)\\
	&= d^\top A^\top Q A d + \frac{b^2}{p^\top Q^{-1} p}\\
	&= d^\top A^\top \left( I - \frac{p (Q^{-1} p)^\top}{p^\top Q^{-1} p} \right) d + \frac{b^2}{p^\top Q^{-1} p}\\
	&= d^\top A d + \frac{b^2}{p^\top Q^{-1} p}\\
	c^\top x^*(d) &= c^\top A d + \frac{b \cdot p^\top Q^{-1} c}{p^\top Q^{-1} p}\\
	\frac{1}{2} x^*(d)^\top Q x^*(d) - c^\top x^*(d) &= \frac{1}{2} d^\top A d - c^\top A d + \frac{b^2/2 - b \cdot p^\top Q^{-1} c}{p^\top Q^{-1} p}.
\end{align*}
Clearly we have $\frac{1}{2} x^*(c)^\top Q x^*(c) - c^\top x^*(d) = - \frac{1}{2} c^\top A c + \frac{b^2/2 - b \cdot p^\top Q^{-1} c}{p^\top Q^{-1} p}$ so therefore
\[ L(d,c) = \frac{1}{2} d^\top A d - c^\top A d + \frac{1}{2} c^\top A c = \frac{1}{2} (d-c)^\top A (d-c). \]
The result now follows.
\Halmos
\endproof

\proof{Proof of Proposition \ref{prop:jpo-portfolio-least-squares}.}
First, notice that $\bbE\left[ \frac{1}{2} \|Vw - c\|_2^2 \right] = \frac{1}{2} \Tr(V^\top V \bbE[w w^\top]) - \Tr(V^\top \bbE[c w^\top]) + \frac{1}{2} \bbE[\|c\|_2^2]$. Via standard vector calculus, we have $\grad_V \bbE\left[ \frac{1}{2} \|Vw - c\|_2^2 \right] = V \bbE[w w^\top] -  \bbE[c w^\top]$, therefore the optimality condition of the least squares predictor is
\[ V \bbE[ww^\top] = \bbE[c w^\top].\]

Now observe that we can write $\bbE\left[ \frac{1}{2} (Vw - c)^\top A (Vw-c) \right] = \frac{1}{2} \Tr(V^\top A V \bbE[w w^\top]) - \Tr(V^\top \bbE[Ac w^\top]) + \frac{1}{2} \bbE[c^\top A c]$. The gradient is $\grad_V \bbE\left[ \frac{1}{2} (Vw - c)^\top A (Vw-c) \right] = A V \bbE[w w^\top] - A\bbE[c w^\top]$. The optimality condition is $A V \bbE[w w^\top] = A \bbE[c w^\top]$. We can alternatively represent this as
\[ V \bbE[w w^\top] \in \bbE[c w^\top] + \left\{ \tilde{V} : A \tilde{V} = \bm{0} \right\}. \]
Since $\bbE[ww^\top]$ is invertible, the result follows.
\Halmos
\endproof

\proof{Proof of Proposition \ref{prop:jpo-knapsack-regularized-mip}}
This result follows immediately from Lemma \ref{lemma:jpo-knapsack-regularized-loss-MIP} below.
\Halmos
\endproof

\begin{lemma}\label{lemma:jpo-knapsack-regularized-loss-MIP}
Assume that $\|d\|_1 \leq M$ and that $M_\tau = M/\left(\min_{j \in [n]} p_j\right)$. The set of optimal solutions to $\max_{x \in X} \left\{ d^\top x - \frac{\lambda}{2} \|x\|_2^2 \right\}$ can be characterized as
\[ \left\{ x : \begin{aligned}
	p^\top x &\leq B, \ \bm{0} \leq x \leq \bm{1}\\
	\tau &\geq 0, \ q,z \in \{0,1\}^n, \ v \in \{0,1\}\\
	\tau &\leq M_\tau v, \ B - p^\top x \leq B (1-v)\\
	d_j - p_j \tau &\leq Mq_j, \ p_j \tau - d_j \leq (M_\tau p_j + M) (1-q_i), \ j \in [m]\\
	d_j - p_j \tau - \lambda &\leq M z_i, \ \lambda + p_j \tau - d_j \leq (M_\tau p_j + M + \lambda) (1-z_j), \ j \in [m]\\
	x_j &\leq q_j, \ x_j \geq z_j, \ j \in [m]\\
	\lambda x_j &\leq d_j - p_j \tau + (M+M_\tau p_j) (1-q_j), \ \lambda x_j \geq d_j - p_j \tau - M z_j, \ j \in [m].
\end{aligned} \right\}. \]
\end{lemma}
\proof{Proof of Lemma \ref{lemma:jpo-knapsack-regularized-loss-MIP}.}
Fix $d$. We consider the primal-dual pair of problems for the regularized fractional knapsack:
\begin{align*}
	&\max_x \left\{ d^\top x - \frac{\lambda}{2} \|x\|_2^2 : p^\top x \leq B, \ \bm{0} \leq x \leq \bm{1} \right\}\\
	&= \min_{s,y,\tau} \left\{ B \tau + \bm{1}^\top y + \frac{1}{2 \lambda} \|s\|_2^2 : s \geq d - p \tau - y, \ s,y,\tau \geq 0 \right\}.
\end{align*}
Using the complementary slackness conditions, the set of primal-dual optimal pairs $(x,s,y,\tau)$ can be written as
\begin{align*}
	X^*(d) &= \left\{ (x,s,y,\tau) : \begin{aligned}
		p^\top x &\leq B, \ \bm{0} \leq x \leq \bm{1},\\
		s &\geq d - p \tau - y, \ s,y,\tau \geq 0,\\
		\tau (B - p^\top x) &= 0,\\
		y_i (1 - x_i) &= 0, \ i \in [n]\\
		x_i (s_i - (d_i - p_i \tau - y_i)) &= 0, \ i \in [n]\\
		s &= \lambda x
	\end{aligned} \right\}\\
	\Proj_{x,\tau}(X^*(d)) &= \left\{ (x,\tau) : \begin{aligned}
		p^\top x &\leq B, \ \bm{0} \leq x \leq \bm{1},\\
		\tau &\geq 0,\\
		\tau (B - p^\top x) &= 0,\\
		\lambda x_i &= \max\left\{0, \min\left\{ \lambda, d_i - p_i \tau \right\} \right\}, \ i \in [n]
	\end{aligned} \right\}.
\end{align*}
To see why the second equality holds, consider some solution $(x,s,y,\tau) \in X^*(d)$. If $d_i - p_i \tau - y_i > 0$, then we need $\lambda x_i = s_i = d_i - p_i \tau - y_i$, which follows from $s_i \geq d_i - p_i \tau - y_i$, $x_i(s_i - (d_i - p_i \tau - y_i))=0$ and $s_i = \lambda x_i$. If $d_i - p_i \tau - y_i \leq 0$, then since $s_i = \lambda x_i$, we would have $x_i(s_i - (d_i - p_i \tau - y_i)) > 0$ if $s_i = \lambda x_i > 0$, so we must have $s_i = \lambda x_i = 0$. Therefore $\lambda x_i = \max\{0,d_i - p_i \tau - y_i\}$. We now show that $y_i = \max\left\{0, d_i - p_i \tau - \lambda \right\}$. To see this, suppose that $d_i - p_i \tau - \lambda > 0$. We know that $y_i \geq d_i - p_i \tau - \lambda x_i > 0$, and since $y_i (1-x_i)=0$, we have $x_i = 1$. Since $\lambda x_i = \lambda = \max\{0,d_i - p_i \tau - y_i\} = d_i - p_i \tau - y_i$ implies that $y_i = d_i - p_i \tau - \lambda$. Now suppose that $d_i - p_i \tau - \lambda \leq 0$. If $y_i > 0$, then since $y_i(1-x_i)=0$, we necessarily have $x_i = 1$. But then $\lambda x_i = \lambda \leq d_i - p_i \tau - y_i$ is a contradiction. Therefore we necessarily have $y_i = 0$. Substituting $y_i = \max\left\{0, d_i - p_i \tau - \lambda \right\}$ into $\lambda x_i = \max\{0,d_i - p_i \tau - y_i\}$ gives us $\lambda x_i = \max\left\{0, \min\left\{ \lambda, d_i - p_i \tau \right\} \right\}$.

We assume that $\|d\|_1 \leq M$, and that we are given an a priori bound $\tau \leq M_\tau$. We can model the constraint $\tau (B - p^\top x) = 0$ as
\[ \tau \leq M_\tau v, \quad B - p^\top x \leq B (1-v), \quad v \in \{0,1\}. \]
We now describe how to model the constraint $\lambda x_i = \max\left\{0, \min\left\{ \lambda, d_i - p_i \tau \right\} \right\}$. First, since $\lambda > 0$, we have that $\min\left\{ \lambda, d_i - p_i \tau \right\} \geq 0$ if and only if $d_i - p_i \tau \geq 0$. Let $q_i$ be an indicator variable for this event, which we model as
\[ d_i - p_i \tau \leq Mq_i, \quad p_i \tau - d_i \leq (M_\tau p_i + M) (1-q_i). \]
Let $z_i$ be an indicator variable for the event $d_i - p_i \tau \geq \lambda$, so we need the constraints
\[ d_i - p_i \tau - \lambda \leq M z_i, \quad \lambda + p_i \tau - d_i \leq (M_\tau p_i + M + \lambda) (1-z_i). \]
Note that implicitly, we have $q_i \geq z_i$. When $q_i = 0$, we have $\lambda x_i = 0$. When $q_i = 1$ and $z_i = 1$, we have $x_i = 1$, and when $q_i=1$ and $z_i = 0$, we have $\lambda x_i = d_i - p_i \tau$. Therefore we need the constraints
\[ x_i \leq q_i, \quad x_i \geq z_i, \quad \lambda x_i \leq d_i - p_i \tau + (M+M_\tau p_i) (1-q_i), \quad \lambda x_i \geq d_i - p_i \tau - M z_i. \]

This shows that the proposed MIP representation is correct.
\Halmos
\endproof

\proof{Proof of Proposition \ref{prop:multiclass-SPOp-loss}.}
The dual of $\min_{x \in X} \left( 2d_j x_j + \sum_{j' \in [m], j' \neq j} (2d_{j'}-1) x_j \right)$ is $\max_{\gamma} \left\{ \gamma : \gamma \leq 2d_j, \ \gamma \leq 2d_{j'} - 1, \ j' \in [m] \setminus \{j\} \right\}$. The lifted representation immediately follows from this.
\Halmos
\endproof

}


\end{document}